\documentclass[a4paper,10pt]{amsart}

\usepackage[english]{babel} \usepackage{amsmath, enumerate, amsfonts,
amssymb, amsthm}

\newtheorem{defin}{\bf Def\mbox{}inition}[subsection]
\newtheorem{algo}[defin]{\bf Algorithm}
\newtheorem{theo}[defin]{\bf Theorem}
\newtheorem{prop}[defin]{\bf Proposition}
\newtheorem{lem}[defin]{\bf Lemma}
\newtheorem{cor}[defin]{\bf Corollary}

\newtheorem{rem}[defin]{\bf Remark}
\newtheorem{ex}[defin]{\bf Example}

\newtheorem*{clai*}{\bf Claim}
\newtheorem*{rem*}{\bf Remark}
\newtheorem*{nota*}{\bf Notation}
\newtheorem*{prop*}{\bf Proposition}

\newcommand{\dps}{\displaystyle}

\def\details#1{#1}

\newcommand{\R}{\mathbb{R}}
\newcommand{\C}{\mathbb{C}}
\newcommand{\N}{\mathbb{N}}
\newcommand{\Z}{\mathbb{Z}}

\renewcommand{\k}{\mathbf{k}}

\newcommand{\pd}{\partial}

\newcommand{\ord}{\mathrm{ord}}
\newcommand{\gr}{\mathrm{gr}}   
\newcommand{\ini}{\mathrm{in}}  
\newcommand{\lm}{\mathrm{lm}}   
\newcommand{\lc}{\mathrm{lc}}   
\newcommand{\lt}{\mathrm{lt}}   
\newcommand{\ND}{\mathrm{Supp}} 
\newcommand{\Exp}{\mathrm{Exp}} 

\newcommand{\E}{\mathcal{E}} 

\newcommand{\RR}{\mathcal{R}}
\newcommand{\CC}{\mathcal{C}}

\newcommand{\Oan}{\mathcal{O}_0}
\newcommand{\Ohat}{{\hat{O}}}
\newcommand{\Oalg}{\k[x]_0} 

\newcommand{\Dan}{\mathcal{D}_0}
\newcommand{\Dhat}{\hat{D}}
\newcommand{\Dalg}{D_{ \{ 0 \} }} 

\newcommand{\uv}{{(u,v)}}
\newcommand{\oo}{{(\mathbf{1,1})}} 
\newcommand{\zo}{{(\mathbf{0,1})}} 
\newcommand{\hooD}{h_\oo(D)}
\newcommand{\hzoD}{h_\zo(D)}
\newcommand{\hzoDalg}{h_\zo(\Dalg)}
\newcommand{\hzoDan}{h_\zo(\Dan)}
\newcommand{\hzoDhat}{h_\zo(\Dhat)}

\newcommand{\Uloc}{\mathcal{U}_{\mathrm{loc}}} 
\newcommand{\Uglob}{\mathcal{U}_{\mathrm{glob}}} 

\newcommand{\Wloc}{\mathcal{W}_{\mathrm{loc}}} 
\newcommand{\Wglob}{\mathcal{W}_{\mathrm{glob}}} 

\newcommand{\supp}{\mathrm{supp}} 
\newcommand{\New}{\mathrm{New}}   
\newcommand{\conv}{\mathrm{conv}} 
\newcommand{\face}{\mathrm{face}}
\newcommand{\Norm}{\mathrm{N}}    

\title[Local Gr\"obner fan]{Local Gr\"obner fan: polyhedral
and computational approach}

\author{Rouchdi Bahloul \and Nobuki Takayama}
\address{Institut Camille Jordan, UMR CNRS 5208,
Universit\'e Claude Bernard - Lyon 1,
B\^at. Jean Braconnier,
43 Boulevard du 11 novembre 1918,
69622 Villeurbanne,
France}
\email{bahloul@math.univ-lyon1.fr}
\address{
Department of Mathematics,
Faculty of Science,
Kobe University,
1-1, Rokkodai, Nada-ku,
Kobe 657-8501, Japan}
\email{takayama@math.kobe-u.ac.jp}

\begin{document}

\begin{abstract}
The goal of this paper is to show that the local Gr\"obner fan 
is a polyhedral fan for ideals in the ring of power series and
the homogenized ring of analytic differential operators. 
We will also discuss about relations between the local Gr\"obner fan and
the (global) Gr\"obner fan for a given ideal
and algorithms of computing local Gr\"obner fans.
In rings of differential operators, finiteness and convexity of local
Gr\"obner fans were firstly proved by Assi, Castro-Jim\'enez and Granger.
But they did not prove that they are polyhedral fans.
\end{abstract}

\maketitle

\tableofcontents

\section*{Introduction}


The Gr\"obner fan for a (homogeneous) polynomial ideal was introduced by
Mora and Robbiano \cite{morarobbiano}. It was also studied by Sturmfels
\cite{sturmfels}. Given an ideal and a weight vector on the variables
we are interested in the equivalence class (the Gr\"obner cone) of the
weight vector defined by the equality of the initial ideals.
In this case it has been proved that the closure of the Gr\"obner cones
form a polyhedral fan. In the ring of algebraic differential operators,
Assi, Castro and Granger \cite{acg00} studied the structure
of the Gr\"obner cones but they did not prove that
it gives rise to a polyhedral fan. This last fact has been proved
by Saito, Sturmfels and the second author \cite{sst}.

Is there a similar construction for power series ideals or more
generally for ideals in the homogenized ring of analytic differential
operators $\hzoDan$?
Assi, Castro-Jim\'enez and Granger \cite{acg01} proved that any standard
cone (Gr\"obner cone) is open, polyhedral and convex and they are in
finite number. Since the ring of analytic power series $\Oan$, as a
subring of $\hzoDan$, is graded, their theorem implies that any
standard cone (Gr\"obner cone) of an ideal in $\Oan$ is open, polyhedral,
and convex. We note that this fact has already been remarked by
Assi \cite{assi}. However, the importance of polyhedral properties of
the fan was not evoked and it will be one of our main purpose here.

The way of discussing about the ring of power series and 
the homogenized ring of analytic differential operators are
analogous. 
However, we will make separate treatments,
because the commutative case is interesting itself and
more readers will be interested in. 

Historically, analogous constructions to Gr\"obner fan already
appeared in works by Lejeune-Jalabert and Teissier \cite{lejeune-t}
with the notion of critical tropism. In $D$-module theory,
we have the notion of slopes with works by Laurent and
Mebkhout \cite{laurent, mebkhout, laurent-meb}, see also \cite{acg96}.
Critical tropisms or slopes are in fact contained in the trace of a
Gr\"obner fan or a standard fan on a $2$-dimensional space.

Now, as an introductive illustration, let us review the case of
a principal ideal in the formal power series ring.

Let $f$ be a formal power series
 $ \sum_{\alpha \in E} c_\alpha x^\alpha $
in $n$ variables.
Here, the set $E$ is the support of the power series
in the space of exponents.
The convex polyhedron
\[\New(f) = \conv \{ \alpha + \N^n  |  \alpha \in E \} \]
is called the Newton polyhedron of $f$.
We are interested in the polyhedral structure of it.

\begin{figure}

\begin{center}
\setlength{\unitlength}{1mm}
\begin{picture}(90,80)
\put(10,40){\circle*{1}} \put(10,36){$-3$}
\put(20,40){\circle*{1}} \put(20,36){$-2$}
\put(30,40){\circle*{1}} \put(30,36){$-1$}

\put(40,40){\circle*{1}} \put(38,41){$0$}  
\put(50,40){\circle*{1}} \put(50,36){$1$}
\put(60,40){\circle*{1}} \put(60,36){$2$}
\put(70,40){\circle*{1}} \put(70,36){$3$}
\put(80,40){\circle*{1}} \put(80,36){$4$}
\put(90,40){\circle*{1}} \put(88,36){$5$}

\put(40,50){\circle*{1}} \put(37,50){$1$}
\put(40,60){\circle*{1}} \put(37,60){$2$}
\put(40,70){\circle*{1}} \put(37,68){$3$}

\put(40,30){\circle*{1}} \put(35,30){$-1$}
\put(40,20){\circle*{1}} \put(35,20){$-2$}
\put(40,10){\circle*{1}} \put(35,10){$-3$}
\put(40,0){\circle*{1}}  \put(35,1){$-4$}

\put(70,40){\line(1,0){20}}
\put(40,60){\line(3,-2){30}}
\put(40,60){\line(0,1){10}}
\put(55,60){${\rm New}(x^3-y^2)$}

\put(40,40){\vector(-1,0){40}}
\put(5,41){{\large $\tau_1$}}

\put(40,40){\vector(0,-1){40}}
\put(42,1){{\large $\tau_3$}}

\put(40,40){\vector(-2,-3){25}}
\put(20,7){$\tau_2$}

\put(20,27){$\sigma_1$}
\put(30,15){$\sigma_2$}

\put(1,17){$\bar{\mathcal{E}}(x^3-y^2)$}

\end{picture}
\end{center}

     \caption{Normal fan and Newton polyhedron for $x^3-y^2$}
     \label{fig:cusp}
\end{figure}
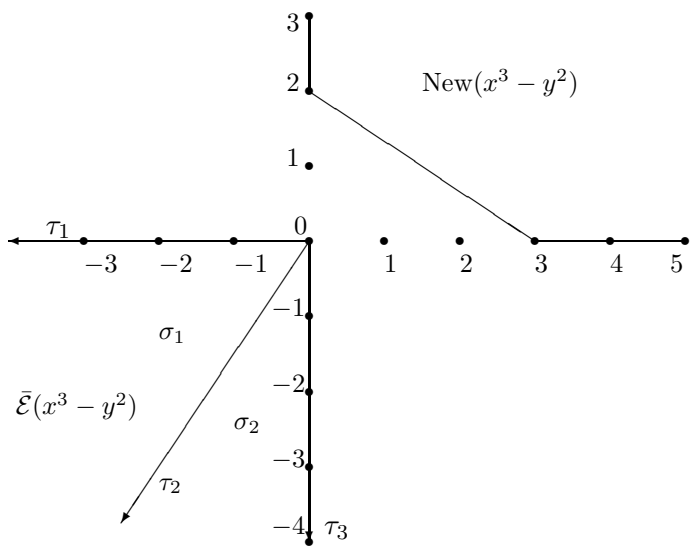

Let $u$ be a vector in $\R_{\leq 0}^n$,
which we will call a weight vector.
For a given weight vector $u$,
let us consider the height of points of the Newton polyhedron
with respect to the direction $u$.
The set of the highest points of $\New(f)$
\[ \face_u (\New(f)) = \{ \alpha \in \New(f) \subset \R^n  | 
u\cdot \alpha \geq u\cdot \alpha' \mbox{ for all } \alpha' \in \New(f) \}
\]
is called the face of $\New(f)$ with respect to $u$.
It is known that there exist only a finite number of faces.
For a given face $F=\face_{u}(\New(f))$, 
the normal cone of $F$ is defined by
\[ \Norm(F) = \{ u' \in \R_{\leq 0}^n  | 
     \face_{u'}(\New(f)) = \face_{u}(\New(f))  \}.
\]
It can be defined in terms of the initial term as follows
\begin{equation}  \label{eq:normalcone}
 \Norm(F) = \{ u' \in \R_{\leq 0}^n  |  \ini_{u'}(f) = {\rm in}_{u}(f), \
   \supp(u') = \supp(u)  \},
\end{equation}
where $\supp(u)=\{i \,| \, u_i \ne 0\}$.
The normal cone $\Norm(F)$ is an open rational convex polyhedral cone.
The collection of the closures of the normal cones 
\[  \bar{\E}(f) =  \{ \overline{\Norm(F)}  | 
       \mbox{ $F$  runs over the faces of $\New(f)$} \}
\]
is called the normal fan of the Newton polyhedron.
Figure \ref{fig:cusp} illustrates the situation for $f(x,y)=x^3-y^2$.
It is well-known and fundamental that
\emph{the normal fan is a polyhedral fan}. Let us recall the
definition.

\begin{defin}  \rm
A polyhedral fan is a finite collection of
(closed) polyhedral cones satisfying the following axioms:
\begin{itemize}
\item[1.]
Every face of a cone is again a cone.
\item[2.]
The intersection of any two cones is a face of both.
\end{itemize}
\end{defin}

Can we generalize the results above to ideals and left ideals
of differential operators?
We can define an analogue of the normal cone
by using an equivalence of initial ideals 
as in (\ref{eq:normalcone}).
Then, it is natural to ask if the collection of the closures of the
cones is a polyhedral fan or not.
The collection has been called the Gr\"obner fan or the standard fan.
The normal fan defined above is the Gr\"obner fan 
of the principal ideal generated by $f$.

In case of non-principal ideals,
we will see that the Gr\"obner fan is locally a common refinement of
the normal fans of a nice set of generators.
Sturmfels introduced this geometric picture of Gr\"obner fan 
in case of homogeneous ideals of the polynomial ring \cite{sturmfels}.
We will apply his idea to local settings with the help of the division
theorem of Assi, Castro-Jim\'enez and Granger \cite{acg01}. This
division takes place in the homogenized ring of analytic (or formal)
differential operators. Such an environment is necessary if we want a
division adapted to a general weight vector. Indeed, a general
division does not exist in the ring $\Dhat$ of formal differential
operators as illustrated with the following example by T. Oaku.

\begin{ex}\label{ex:impossible-division} \rm
Take $n=1$. We want to divide $f=x$ by $g= x+x^k \pd$ with respect to
a weight vector $w=(u,v)$, with $u<0$ and $u+v\ge 0$. Take $k \ge 2$
sufficiently large such that $x \succ x^k \pd$. Here $\prec$
is a refinement of the partial order defined by $w$. Suppose that
a division theorem holds:
\[  x = q \cdot (x + x^k \pd) + r, \quad q,r \in \Dhat \]
where $r$ cannot be divided by $g$, i.e. the leading monomial of $r$
is not divisible by $x$ the leading monomial of $g$.
But $r=x- q \cdot (x + \pd x^k + k x^{k-1}) =
(1-q - q \pd x^{k-1} -k x^{k-2}) \cdot x$ therefore the leading
monomial of $r$ contains $x$.
It is a contradiction.
\end{ex}

Finally, let us mention about applications of our result.
The first author gave a constructive method to compute
a local Bernstein-Sato polynomial for a given set of analytic
functions or polynomials in \cite{Compos}.
An essential step of the construction is 
the computation of an analytic standard fan,
and we will give algorithms to compute it.
The second application is given by N. Touda who constructs
{\it local tropical variety} based on our result that
``the local Gr\"obner fan is a polyhedral fan'' \cite{touda}.
This is an analogous construction of {\it tropical variety}
by D. Speyer and B. Sturmfels \cite{ss-trop}.

\section{Definitions and statement of the main results}

We fix an integer $n \ge 1$. We denote by $x$ the system of variables
$(x_1,\ldots,x_n)$. The partial derivation $\frac{\partial}{\partial x_i}$
shall be denoted by $\pd_i$. We will use the notation $x^\alpha$ for $x_1^{\alpha_1} \cdots
x_n^{\alpha_n}$ where $\alpha\in \N^n$. In the same way $\pd^\beta$
denotes $\pd_1^{\beta_1} \cdots \pd_n^{\beta_n}$.

The symbol $\k$ denotes a field.
$\Oalg$ denotes the localization of $\k[x]$ at the point $0\in \k^n$.
When $\k=\C$ this is a subring of $\Oan$ germ at $0$ of analytic functions.
The latter is itself a subring of the formal power series
ring $\Ohat=\k[[x]]$.
In this section, the symbol $\CC$ shall denote one of the rings
$\k[x]$, $\Oalg$, $\Oan$, $\Ohat$.

Now let us denote respectively by $D$, $\Dalg$, $\Dan$ and $\Dhat$
the ring of differential operators with coefficients
in $\k[x]$, $\Oalg$, $\Oan$ and $\Ohat$. The symbol $\RR$ shall denote
one of these rings.

Recall that these rings are characterized by the following relations:
for any $c(x) \in \CC$ and $i=1, \ldots, n$, $\pd_i c(x) - c(x) \pd_i
= \frac{\pd c(x)}{\pd x_i}$.

We denote by $\Uloc$ the set $\{u \in \R^n  |  u_i\le 0, \forall i\}$.
This is a subset of $\Uglob=\R^n$. We denote by $\Wloc=\{w=(u,v) \in
\R^{2n}  |  u_i \le 0, u_i+v_i \ge 0, \forall i\}$. This is subset of
$\Wglob= \{w=(u,v) \in\R^{2n} |u_i+v_i \ge 0, \forall i\}$. For each $i$,
$u_i$ shall be seen as a weight on the variable $x_i$ while $v_i$ shall
be a weight on $\pd_i$. The use of the subscripts ``loc'' and ``glob''
together with the relations $u_i \le 0$ and $u_i+v_i \ge 0$ shall be
made clearer in Remark \ref{rem:explanation} (see below).

For the construction of Gr\"obner fans, we will need a division process
with respect to ``arbitrary'' weight vectors in rings of differential
operators. We have seen in Example \ref{ex:impossible-division} that
a general division does not exist e.g. in $\Dhat$.
Therefore, homogenization techniques will be necessary. They
are differential generalizations of that by Lazard \cite{lazard}.

First let us define the $\zo$-homogenization.
Let us introduce the ring $\hzoDhat$. It is the $\k$-algebra generated
by $\Ohat$, $\pd_i$ ($i=1,\ldots,n$) and a new variable $h$ with
\[\pd_i c(x) -c(x) \pd_i = \frac{\pd c(x)}{\pd x_i} \cdot h,\]
$i=1,\ldots,n$, as the only non trivial commutation relations.
We can replace $\Ohat$ by $\k[x], \Oalg, \Oan$ and we obtain the
subrings $\hzoD \subset \hzoDalg \subset \hzoDan \subset \hzoDhat$.
An element $P$ in one of these rings has a unique writing as
$P=\sum_{\alpha, \beta, k} c_{\alpha, \beta, k} x^\alpha \pd^\beta h^k$.
We define the support $\ND(P)\subset \N^{2n+1}$ of $P$ as the set
of $(\alpha, \beta, k)$ such that $c_{\alpha,\beta,k} \ne 0$.
We define its degree $\deg(P)$ as the maximum of $|\beta|+k$.
We say that $P$ is homogeneous if for any $(\alpha, \beta,k)$ in
$\ND(P)$, $|\beta|+k=\deg(P)$. An ideal in $h_\zo(\RR)$ is said
to be homogeneous if it can be generated by homogeneous elements.

Let $P=\sum_{\alpha, \beta} c_{\alpha, \beta} x^\alpha \pd^\beta$ in $\RR$.
The support $\ND(P) \subset \N^{2n}$ is defined in the same way.
We define the $\zo$-degree of $P$ as the maximum of $|\beta|$
for $(\alpha,\beta) \in \ND(P)$. Its $\zo$-homogenization is defined by
$h_\zo(P)=\sum c_{\alpha, \beta} x^\alpha \pd^\beta h^{\deg(P)-|\beta|}$.
It is a homogeneous element of $h_\zo(\RR)$ of degree the degree of $P$.

Now let us introduce the $\oo$-homogenized Weyl algebra $\hooD$
(\cite{kan}, \cite{castro-narvaez}).
It is the $\k$-algebra generated by the symbols $x_i$, $\pd_i$ and $h$ where
the only non trivial relations are
\[\pd_i x_i -x_i \pd_i = h^2; \ i=1, \ldots,n.\]
For $P\in \hooD$, its degree shall be the total degree in all the variables.
It is homogeneous if all its monomials have the same degree.
For $P=\sum_{\alpha, \beta} c_{\alpha, \beta} x^\alpha \pd^\beta$ in $D$,
the $\oo$-degree is the total degree in the variables $x_i$ and $\pd_i$.
We define the $\oo$-homogenization of $P$ as $h_\oo(P)=\sum c_{\alpha,
\beta} x^\alpha \pd^\beta h^{d-|\alpha|-\beta}$ where $d$ is the
$\oo$-degree of $P$. This is a homogeneous element of degree $d$.

For $f\in \CC$ and a local weight vector $u \in \Uloc$,
we denote by $\ord^u(f)$ the maximum of
the scalar products $u \cdot \alpha$ for $\alpha \in \ND(f)$. In the case
$f\in \k[x]$, this definition is valid for $u\in \Uglob$. The order
gives rise to a filtration $F^u(\CC)$ given by $F^u_k=\{f | \ord^u(f)
\le k\}$ as well as to the associated graded ring $\gr^u(\CC)=\oplus_k
F_k^u/ F^u_{<k}$. If $\CC=\k[x]$ and $u\in \Uglob$, this graded ring
is equal to $\k[x]$. All these rings can be seen as subrings of
$\gr^u(\Ohat)$. For $f\in \CC$, its initial form $\ini_u(f)$ shall be
the class of $f$ in $F_k^u/ F^u_{<k}$ where $k=\ord^u(f)$.
In literature, we also find the terminology principal symbol
with the notation $\sigma^u(f)$.
Now if $I$ is an ideal in $\CC$, the restriction of $F^u(\CC)$ to $I$
gives rise to a graded ideal in $\gr^u(\CC)$ which we denote by
$\ini_u(I)$. This ideal is in fact generated by all the $\ini_u(f)$ for
$f\in I$.

Now all the definitions above are valid if we replace $\CC$ by $\RR$ and
$\Uloc$ (resp. $\Uglob$) by $\Wloc$ (resp. $\Wglob$). For example if write
$w=(u,v)$ then $\ord^w(P)$ shall be the maximum of $(u,v) \cdot (\alpha,
\beta)$ for $(\alpha,\beta)\in \ND(P)$. Moreover if we put the weight $0$
to the variable $h$, we can extend the definitions above to
$\hooD$ or $\hzoD$ with $\Wglob$ and to $\hzoDalg$, $\hzoDan$,
$\hzoDhat$ with $\Wloc$.

\begin{rem}\label{rem:explanation}
The subscript ``loc'' means that we are working with local rings. In
this case, the weights $u_i$ of the $x_i$ variables need to be $\le 0$
if we want the definition above to work.
Moreover the conditions $u_i+v_i \ge 0$ are necessary to respect the
non commutative structure of the rings
in the differential case.
Indeed we have: $\ord^w(x_i \pd_i)= u_i+v_i$ and 
$\ord^w(x_i) +\ord^w(\pd_i)=u_i + v_i$. Similarly, we would like to
have $\ord^w(\pd_i x_i)=v_i +u_i$. Since $\pd_i x_i= x_i \pd_i +1$
and $\ord^w(x_i \pd_i)=u_i+v_i$ and $\ord^w(1)=0$, we need $u_i+v_i
\ge 0$.
\end{rem}

Let $u$ be in $\Uloc$ and define a strata in the weight space:
\[S_\CC(u)=\{ u \in \Uloc | \gr^u(\CC)=\gr^{u'}(\CC)\}.\]
We simply write $S(u)$ when no confusion arises.
Then for a given ideal $I$ in $\CC$, consider the equivalence relation
\[u \sim u' \iff u'\in S(u) \textrm{ and } \ini_u(I)=\ini_{u'}(I)\]
for which we denote by $C_I[u]$ the class of $u$. This class is called
a Gr\"obner cone (the Gr\"obner cone of $I$ w.r.t. $u$).
We then denote by $\E(I,\Uloc)$ the set the Gr\"obner cones $C_I[u]$
for $u\in \Uloc$. This is the \emph{open Gr\"obner fan} of $I$ on $\Uloc$.
If $\RR=\k[x]$, we can do the same by replacing $\Uloc$ with $\Uglob$ and
obtain $\E(I,\Uglob)$. When the context is clear we shall only write
$\E(I)$ (for example if $\CC=\Ohat$).
The notation $\E$ comes from the word ``\'eventail''.

Finally, the closure of $C_I[u]$ will be called the closed Gr\"obner
cone of $I$ w.r.t. $u$. The set of the closed Gr\"obner cones
shall be denoted by $\bar{\E}(I,\Uloc)$ (resp. $\bar{\E}(I,\Uglob)$)
and called \emph{the closed Gr\"obner fan} of $I$.

We can give the same definitions for an ideal $I$ in $\RR$, $\hooD$ or
$h_\zo(\RR)$.\\

\noindent
{\bf Note.}
\begin{itemize}
\item
In \cite{morarobbiano} Mora and Robbiano introduced and studied
the Gr\"obner fan for a homogeneous ideal in $\k[x]$. It has been
proved by Sturmfels \cite{sturmfels} that the closed Gr\"obner fan
is a polyhedral fan.
\item
In \cite{acg00}, Assi, Castro-Jim\'enez and Granger introduced $\E(h_\oo(I),
\Wglob)$ for an ideal $I$ in $D$ and called it the Gr\"obner fan of $I$.
It is not strictly speaking a fan, although for this case, Saito, Sturmfels
and Takayama \cite{sst} proved that $\bar{\E}(J)$ is a polyhedral fan
for any homogeneous $J$ in $\hooD$.
\item
In \cite{acg01}, Assi, Castro-Jim\'enez and Granger defined and worked with
$\E(h_\zo(I), \Wloc)$ for $I$ in $\Dan$ (or $\Dhat$). They proved it is
a finite collection of convex rational polyhedral cones and called it the
(analytic) standard fan of $h_\zo(I)$, however they did not prove that
$\bar{\E}(h_\zo(I))$ is a polyhedral fan.
\end{itemize}

Let us state the main contributions of the present paper.
\begin{itemize}
\item[(A)]
Let $I$ be an ideal in $\Oan$ or $\Ohat$ then the closed Gr\"obner fan
$\bar{\E}(I, \Uloc)$ is a polyhedral fan.
\item[(B)]
Let $I$ be a homogeneous ideal in $\hzoDan$ or $\hzoDhat$ then the closed
Gr\"obner fan $\bar{\E}(I,\Wloc)$ is a polyhedral fan.
\item[(C)]
Let $I$ be an ideal in $k[x]$ then $\bar{\E}(\Oalg I, \Uloc) =
\bar{\E}(\Ohat I, \Uloc)$ and we provide an algorithm to compute its
restriction to any linear subspace.
\item[(D)]
Let $I$ be an ideal in $D$ then we provide an algorithm to compute
the open Gr\"obner fan of $\Dhat I$ and the closed Gr\"obner fan of 
$h_\zo(\Dhat I)$ (the latter being a fan).
\item[(E)] Comparison theorems of local and global Gr\"obner fans 
(Theorems  \ref{cor:Wloc'}, \ref{prop:locfan}, \ref{cor:Wloc'Hom},
\ref{prop:locfanH}, \ref{theo:glob=locIFhomo}).
\end{itemize}

The fan in (C) (resp. (D)) shall be called the local
Gr\"obner fan (at $x=0$) of $I$.

The way of discussing about the ring of power series and 
the homogenized ring of analytic differential operators are
analogous. 
However, we will make separate treatments;
in section 2, we discuss about the ring of power series
and, in section 3, we only
prove the results which are specific to the case of differential
operators and we refer to section 2 otherwise.
In section 4 we give general results about the different kinds of
Gr\"obner fans and we compare them (E). In section 5 we will be concerned
with (C) and (D). Results (A) and (B) will be fundamental for enumerating
the Gr\"obner cones. We end section 5 with some examples.

\section{The closed Gr\"obner fan in $\Ohat$ or $\Oan$ is a polyhedral fan}

In this section we will prove the following theorem.

\begin{theo}\label{theoC:SCGCpolyhedral}
Let $I$ be an ideal in $\Oan$ or in $\Ohat$ then the closed Gr\"obner
fan $\bar{\E}(I,\Uloc)$ is a rational polyhedral fan.
\end{theo}

\subsection{Division theorem and standard bases in $\gr^u(\Ohat)$
and $\gr^u(\Oan)$}

In the following, we will state a division theorem in the graded
rings $\gr^u(\Ohat)$ and $\gr^u(\Oan)$, let us describe them
more explicitly.

\begin{lem}\label{lemC:grCC}
Let $u \in \Uloc$ be a local weight vector and
for some $0\le m\le n$ assume that $u=(u_1,\ldots,u_m, 0, \ldots ,0)$
with $u_i\ne 0$. Then\\
$\gr^u(\Oan)=
\C\{x_{m+1}, \ldots, x_n\}[x_1,\ldots,x_m]$ and
$\gr^u(\Ohat)=\k[[x_{m+1}, \ldots, x_n]][x_1, \ldots,x_m]$.
\end{lem}

The proof immediately follows from the definition. 
Let us denote by $\CC$ one
of $\Ohat$, $\Oan$. Thanks to the preceding lemma, we shall see
$\gr^u(\CC)$ as a subring of $\CC$. An element $f \in \CC$ is said
to be $u$-homogeneous (or simply homogeneous if the context is clear)
if all its monomials have the same $u$-order.

Let $\prec$ be a total order on $\N^n$ (or equivalently on the terms
$x^\alpha=x_1^{\alpha_1} \cdots x_n^{\alpha_n}$). It is said to be a
monomial order if $\alpha \prec \alpha'$ implies $\delta + \alpha \prec
\delta + \alpha'$. It is called local if moreover $\alpha \preceq 0$ for
any $\alpha$. On the opposite side, a monomial order with
$\alpha \succeq 0$ is called a global or a well order.

Let $f\in \CC$ be a non zero power series then we define
its leading exponent $\exp_\prec(f)$ w.r.t. $\prec$ as the maximum
of $\ND(P)$ for $\prec$. We define its leading term $\lt_\prec(f)=
x^{\exp_\prec(f)}$, its leading coefficient $\lc_\prec(f)$
(in $\k$ or $\C$) as the coefficient corresponding to $\lt_\prec(f)$ and
its leading monomial $\lm_\prec(f)=\lc_\prec(f) \lt_\prec(f)$.

Let us state a division theorem in $\gr^u(\CC)$ for $\prec$ with unique
quotients and remainder similarly to \cite{cg}.

Let $g_1,\ldots,g_r$ be non zero homogeneous elements in $\gr^u(\CC)$.
Define a partition of $\N^n$ associated with the $\exp_\prec(g_j)$ as
$\Delta_1=\exp_\prec(g_1)+ \N^n$,
$\Delta_j=(\exp_\prec(g_j) +\N^n) \smallsetminus (\Delta_1 \cup
\cdots \cup \Delta_{j-1})$ for $j>1$ and
$\bar{\Delta}=\N^n \smallsetminus  (\cup_j \Delta_j)$.

\begin{theo}[Division theorem in $\gr^u(\Oan)$ and $\gr^u(\Ohat)$]
\label{theo:div}
For any homogeneous $f\in \gr^u(\CC)$, there exists a unique
$(q_1,\ldots, q_r, R) \in \gr^u(\CC)^{r+1}$ made of homogeneous elements
such that $f=\sum_j q_j g_j + R$ and
\begin{itemize}
\item
for any $j$, either $q_j=0$ or $\ND(q_j) +\exp_\prec(g_j) \subset \Delta_j$,
\item
either $R=0$ or $\ND(R) \subset \bar{\Delta}$.
\end{itemize}
We call $R$ the remainder of the division of $f$ by the $g_j$ w.r.t. $\prec$.
\end{theo}

\begin{proof}
Let us first suppose $u=(0)$. In this case $\gr^u(\CC)=\CC$ and
any element in $\CC$ is $u$-homogeneous (of order $0$).
Let us show how we can derive the result from \cite[Theorem 1.5.1]{cg}.
By Robbiano's theorem \cite{robbiano}, the order $\prec$ can be defined
as a lexicographical order with respect to $n$ weight vectors.
Let $w^1$ be the first of them. Then define the order $<_{w^1}$ in a
lexicographical way by $w^1$ and by $<_0$ where $<_0$ is the inverse
of $\prec$ (and thus $<_0$ is a well order). The orders $\prec$ and
$<_{w^1}$ are equivalent. Since $\prec$ is local $w^1$ has non positive
components. Thus with the order $<_{w^1}$, we are exactly under the
hypotheses of \cite[Theorem 1.5.1]{cg}. Now suppose $u$ is not
necessarily zero. Since $\gr^u(\CC)$ is a subring of $\CC$ (by lemma
\ref{lemC:grCC}), we make
the division in $\CC$ and we remark that the division process
conserves $u$-homogeneity which concludes the proof.
\end{proof}
\noindent
As a consequence we have:

$\exp_\prec(f)=\max_\prec\{ \exp_\prec(q_j g_j), j=1,\ldots,r;
\exp_\prec(R)\}$.\\
With the notations of the proof above, this implies the following:

$\ord^{w^1}(f)=\max\{\ord^{w^1}(q_j g_j), j=1,\ldots,r; \ord^{w^1}(R) \}$.\\

Let $J$ be an ideal in $\gr^u(\CC)$. We suppose $J$ to be $u$-homogeneous,
i.e. generated by homogeneous elements. We define the set of the leading
exponents of $J$ as
\[\Exp_\prec(J)=\{ \exp_\prec(f) | f \in J, \ f\ne 0\}.\]
This set is stable by sums in $\N^n$ thus by Dickson lemma:

\begin{defin} \rm
There exists $G=\{g_1,\ldots, g_r\} \subset J$ (made of homogeneous
elements) such that $\Exp_\prec(J)= \cup_j (\exp_\prec(g_j) + \N^n)$.
Such a set $G$ is called a (homogeneous) $\prec$-standard basis of $J$.
\end{defin}

The statement concerning homogeneity follows from the following:
take $g_j\in G$, then any $u$-homogeneous part of $g$ belongs to $J$.
Let $g'_j$ be the one that contains the leading term of $g$.
The set of $g'_j$ is a homogeneous standard basis.

The following statements are equivalent, see \cite[Cor. 1.5.4]{cg}:
\begin{itemize}
\item
A set $G\subset J$ is a homogeneous standard basis of $J$
\item
For any homogeneous $f\in \gr^u(\CC)$: $f\in J$ if and only if the
remainder of the division of $f$ by $G$ is zero.
\end{itemize}
Keeping the notations of the proof above, we see that for such a standard
basis, the set of $\ini_{w^1}(g)$ for $g\in G$ generates $\ini_{w^1}(J)$.
In fact we have more than that: Given a local order $\prec$ and $u\in
\Uloc$, we define the local order $\prec_u$ in a lexicographical way
by $u$ and $\prec$.

\begin{lem}\label{lemC:SB->in}
\begin{itemize}
\item[(1)]
Given $I \subset \CC$, if $G$ is a $\prec_u$-standard basis of $I$
then $\ini_u(G)$ is a (homogeneous) standard basis of $\ini_u(I)$
for $\prec$.
\item[(2)]
$\Exp_{\prec_u}(I)=\Exp_\prec(\ini_u(I))$.
\end{itemize}
\end{lem}

\begin{proof}
Statement (2) easily follows from (1) and the fact that for any
$g\in G$, $\exp_{\prec_u}(g)= \exp_\prec(\ini_u(g))$. Let us prove (1).
Let $f'$ be in $\ini_u(I)$. We want to prove that $\exp_\prec(f')$
is in $\exp_\prec(\ini_u(g))+\N^n$ for some $g\in G$.
By considering
the homogeneous part of $f'$ that contains the $\prec$-leading
term, we may assume $f'$ to be homogeneous.
Let $f\in I$ be such that $\ini_u(f)=f'$ and let us consider the
division of $f$ by $G$ w.r.t. $\prec_u$ as in Theorem \ref{theo:div}:
$f=\sum_1^r q_j g_j$ where $G=\{g_1, \ldots, g_r\}$ with
\[(\star) \qquad \ND(q_j) +\exp_{\prec_u}(g_j) \subset \Delta_j.\]
We have $\ord^u(f)\ge \ord^u(q_j g_j)$ for any $j$. Let
$q'_j=\ini_u(q_j)$ if $\ord^u(f)= \ord^u(q_j g_j)$ and $q'_j=0$
otherwise then $f'=\sum_1^r q'_j \ini_u(g_j)$.
Since $\exp_{\prec_u}(g)=\exp_\prec(\ini_u(g))$, the partition of
$\N^n$ associated with the $\exp_{\prec_u}(g_j)$ is the same as that
associated with $\exp_\prec(\ini_u(g_j))$. Moreover, for any $h\in
\CC$, $\ND(\ini_u(h)) \subset \ND(h)$, thus the relation $(\star)$
becomes
\[\ND(q'_j) +\exp_{\prec}(\ini_u(g_j)) \subset \Delta_j.\]
This means that the remainder of the division of $f'$ by the
$\ini_u(g_j)$ w.r.t. $\prec$ is zero from which the conclusion
follows.
\end{proof}

Now given a standard basis $G$ of $J \subset \gr^u(\CC)$. We say that $G$
is minimal if for $g,g' \in G$, $\exp_\prec(g) \in \exp_\prec(g')+\N^n$
implies $g=g'$. We say that $G$ is {\bf reduced} if it is minimal,
unitary (i.e. $\lc_\prec(g)=1$ for any $g\in G$) and if for any
$g\in G$, $\ND(g) \smallsetminus \{\exp_\prec(g)\} \subset \N^n
\smallsetminus \Exp_\prec(I)$.

\begin{lem}\label{lemC:reducedSB}
\begin{itemize}
\item[(1)]
Given a homogeneous ideal $J$ in $\gr^u(\CC)$ and a local order $\prec$,
a reduced standard basis w.r.t. $\prec$ exists and is unique. Moreover
it is made of homogeneous elements.
\item[(2)]
Given $I$ in $\CC$, if $G$ is the reduced $\prec_u$-standard basis
then $\ini_u(G)$ is the reduced $\prec$-standard basis of $\ini_u(I)$.
\end{itemize}
\end{lem}

\begin{proof}
The proof of the first statement is classical, we omit it. The second
one easily follows from the fact that for any $g \in G$, $\exp_{\prec_u}(g)=
\exp_\prec(\ini_u(g))$ and $\ND(\ini_u(g)) \subset \ND(g)$.
\end{proof}

\begin{lem}\label{lemC:prec_12}
Suppose we have two local orders $\prec_1$ and $\prec_2$ and a homogeneous
ideal $J\subset \gr^u(\CC)$. Let $G$ be a (resp. the reduced)
$\prec_1$-standard basis of $J$ and suppose for any $g\in G$ that
$\exp_{\prec_1}(g)= \exp_{\prec_2}(g)$ then $G$ is a (resp. the
reduced) $\prec_2$-standard basis of $J$.
\end{lem}

\begin{proof}[Sketch of Proof]
Let $f \in J$. Divide $f$ by $G$ w.r.t $\prec_1$ for which the remainder
is zero. By assumption, this division is also the division w.r.t
$\prec_2$, and since the remainder is zero, $G$ is a $\prec_2$-standard
basis of $J$. The statement concerning the reducibility is trivial.
\end{proof}

\subsection{Back to Gr\"obner fans}

For a local weight vector $u \in \Uloc$,
call $\supp(u)=\{1\le i\le n | u_i \ne 0\}$ \emph{the support of $u$}.
By lemma \ref{lemC:grCC}, we have:
\[u'\in S_\CC(u) \iff \supp(u)=\supp(u') .\]

\begin{prop}\label{propC:GrobnerCone_SB}
Given an ideal $I$ in $\CC$ and $u\in \Uloc$, then for any local order
$\prec$, the $\prec_u$-reduced standard basis $G$ of $I$ satisfies:
\[C_I[u]=\{ u'\in \Uloc | \supp(u)=\supp(u') \textrm{ and }
\forall g\in G, \ini_u(g)=\ini_{u'}(g)\}.\]
\end{prop}

\begin{proof}
We already know that $u'\in S_\CC(u) \iff \supp(u)=\supp(u')$ then under
this assumption we have to prove that $\ini_u(g)=\ini_{u'}(g)$ for any
$g \in G$ if and only if $\ini_u(I)=\ini_{u'}(I)$.
Take $u'$ in the RHS. By assumption, for any $g\in G$, $\exp_{\prec_u}(g)=
\exp_\prec(\ini_u(g))=\exp_\prec(\ini_{u'}(g))= \exp_{\prec_{u'}}(g)$
so by lemma \ref{lemC:prec_12}, $G$ is also the $\prec_{u'}$-reduced
standard basis of $I$. It follows that $\ini_u(I)$ and $\ini_{u'}(I)$
have a common set of generators which is $\ini_u(G)=\ini_{u'}(G)$.

Now take $u'$ in the LHS. By Lemma \ref{lemC:SB->in}(2), we have
$\Exp_{\prec_u}(I)=\Exp_\prec(\ini_u(I))= \Exp_\prec(\ini_{u'}(I))=
\Exp_{\prec_{u'}}(I)$ which implies that $G$ is the reduced standard
basis for $\prec_{u'}$. Thus by lemma \ref{lemC:reducedSB}(2), $\ini_u(G)$
and $\ini_{u'}(G)$ are both the $\prec$-reduced standard basis of
$\ini_u(I)=\ini_{u'}(I)$. By unicity of the reduced standard basis,
$\ini_u(G)=\ini_{u'}(G)$. As a consequence, $u'\in$ RHS.
\end{proof}

\begin{cor}\label{corC:FormAnalGF}
Let $I$ be an ideal in $\Oan$ and let $\hat{I}=\Ohat \cdot I$ then the
Gr\"obner fans $\E(I,\Uloc)$ and $\E(\hat{I},\Uloc)$ are equal.
\end{cor}
\begin{proof}
Let $u$ be in $\Uloc$. Let $G$ be the $\prec_u$-reduced standard
basis of $I$ then by using the Buchberger criterion involving the
$S$-functions (see \cite[Prop. 1.6.2]{cg}), we obtain that $G$ is a
$\prec_u$-standard basis of $\hat{I}$ which in turn implies that it is
the reduced one. Moreover by lemma \ref{lemC:grCC}, we have that
$S_{\Ohat}(u)= S_{\Oan}(u)$ thus $C_I[u]=C_{\hat{I}}[u]$ which completes
the proof.
\end{proof}

With this corollary, we may assume that our ideal is in $\Ohat$. We will
work in $\Ohat$ for the rest of this section.
Let $g$ be in $\Ohat$, we define its Newton polyhedron as the following
convex hull:
\[\New(g)=\conv(\ND(g)+ \N^n).\]
The set $\ND(g)+\N^n$ is by definition stable by sums so let $E(g)$ be
the minimal finite subset of $\ND(g)$ such that
$\ND(g)+\N^n=E(g)+\N^n$.
As a consequence we have
\[\New(g)=\conv(E(g))+\R_{\ge 0}^n\]
which implies that $\New(g)$ is strictly speaking
a polyhedron \cite[p. 30]{ziegler}.

First we have the following result which is a part of
Theorem \ref{theoC:SCGCpolyhedral}.

\begin{theo}
For any ideal $I$ in $\Ohat$, the set $\{\ini_u(I)| u\in \Uloc\}$ is finite.
\end{theo}
\begin{proof}
We sketch the proof by giving the main steps since it is very similar
to that of \cite[Th. 4]{acg01}.
\begin{itemize}
\item[(a)]
Given $g \in \Ohat$, the set of $\exp_\prec(g)$, $\prec$ being any local
order, is finite. By using the same arguments as in the proof of
\cite[Prop. 17]{acg01}, we can easily show that it is contained in $E(g)$.
\item[(a')] The set of $\ini_u(g)$, $u\in \Uloc$, is finite. Indeed it
is in one to one correspondence with the set of faces of $\New(g)$.
\item[(b)]
Given an ideal $I$ in $\Ohat$, the set $\{\Exp_\prec(I) | \prec
\textrm{ is a local order}\}$ is finite. The proof uses the same
arguments as that of \cite[Th. 5]{acg01} and is based on (a). We can
also prove it as \cite[Th. 1.2]{sturmfels}.
\item[(c)]
This last step is the main part.
Thanks to the preceding statement, it is enough
to prove the following: Given $E=\Exp_\prec(I)$ for some local order
$\prec$, the set of all $\ini_u(I)$ such that $\Exp_{\prec_u}(I)=E$
is finite. Let $G$ be a $\prec$-standard basis of $I$ then
for any $u$ as above, $\ini_u(I)$ is generated by the set $\{\ini_u(g)|
g\in G\}$. By (a'), there is only a finite number of such sets.
\end{itemize}
\end{proof}

Now let $u$ be in $\Uloc$ and let $\prec$ be any local order. Let $G$ be
the $\prec_u$-reduced standard basis of $I$. Let
\[Q=\sum_{g\in G} \New(g)\]
be the Minkowski sum of the Newton polyhedra of the $g$ in $G$.

\begin{prop}\label{propC:GrobnerCone_NC}
\[C_I[u]=\Norm_Q(\face_u(Q))\]
where $\Norm_Q(\face_u( \cdot ))$ denotes the normal cone of the face
w.r.t. $u$.
\end{prop}
As a direct consequence:

\begin{cor}
The Gr\"obner cone $C_I[u]$ is a convex rational relatively open
polyhedral cone.
\end{cor}
\noindent
The proof of the proposition will be decomposed into several lemmas.

For $g \in \Ohat$ and $u\in \Uloc$, let $E_u(g)$ be the elements
of $E(g)$ which are maximum for the scalar product with $u$. Let $e_i \in
\N^n$, $i=1,\ldots,n$ be the canonical base of the semi-group $\N^n$.

\begin{lem}\label{lemC:faceNew}
\[\face_u(\New(g))=\conv(E_u(g))+ \sum_{i\notin \supp(u)}
\R_{\ge 0} \cdot e_i.\]
\end{lem}

\begin{proof}
It is easy to see that $\face_u(\R_{\ge 0}^n)$ is equal to the second member
of the RHS. Moreover there is a general identity for polyhedra
\[
\face_u(P+P')=\face_u(P)+ \face_u(P').
\]
Thus it suffices to show this equality: $\face_u(\conv(E(g)))=\conv(E_u(g))$.
Let us prove it.
Take $\alpha$ in the LHS. Written as a sum $\alpha=\sum_j c_j \alpha_j$,
$c_j\ge 0$, $\sum_j c_j=1$, of elements in $E(g)$, we see that the
scalar product $u\cdot \alpha$ is equal to some $u \cdot \alpha_j$.
This implies that $u$ evaluated on the LHS or on the RHS gives the
same number. From this observation, we easily derive the desired
equality.
\end{proof}

Let $g_1,\ldots,g_r$ be in $\Ohat$. Put $Q_j=\New(g_j)$ and let $Q$ be the
Minkowski sum of the $Q_j$'s.

\begin{lem}\label{lemC:NormalFace}
For any $u\in \Uloc$, we have
\[\Norm_Q(\face_u(Q))= \bigcap_{j=1}^r \Norm_{Q_j}(\face_u(Q_j)).\]
\end{lem}

This lemma is known for polytopes (see \cite{ziegler}, Prop. 7.12, p. 198).
We can prove it by reducing the case of polyhedra to the case of polytopes
by truncating our polyhedra.

Now, here is the last lemma before the proof of proposition
\ref{propC:GrobnerCone_NC}.

\begin{lem}\label{lemC:ini=face}
Let $g\in \Ohat$ then for $u,u'\in \Uloc$:
\[\big[ \ini_u(g)=\ini_{u'}(g) \textrm{ and } \supp(u)=\supp(u') \big] \iff
\face_u(\New(g))= \face_{u'}(\New(g)).\]
\end{lem}

\begin{proof}
First let us prove this claim:\\
{\bf Claim.} Under the condition $\supp(u)=\supp(u')$ the following
equivalence holds: $E_u(g)=E_{u'}(g) \iff \ini_u(g)=\ini_{u'}(g)$.\\
Suppose we have $E_u(g)=E_{u'}(g)$,
it is enough to prove $\ND(\ini_u(g))= \ND(\ini_{u'}(g))$.
Let $\alpha$ be in $\ND(\ini_u(g))$. Then there exists $e\in E(g)$
and $\alpha' \in \N^n$ such that $\alpha= e+\alpha'$.
We have $u\cdot \alpha \le u\cdot e$ but by definition of
$\ini_u(g)$, $u\cdot \alpha\ge u\cdot e$. Thus $e$ belongs to $E_u$ which
is equal to $E_{u'}(g)$ by assumption. This also means that $u \cdot \alpha'=0$
but since $\supp(u)=\supp(u')$, this implies $u' \cdot \alpha'=0$.
Thus $u' \cdot \alpha=u' \cdot e$. Now take any $\alpha''$ in $\ND(g)$.
We can write $\alpha''=e'+a'$ with $e'\in E$ and $a'\in \N^n$.
Therefore, $u' \cdot \alpha'' \le u' \cdot e' \le u'\cdot e=u'\cdot \alpha$.
This means that $\alpha$ belongs to
$\ND(\ini_{u'}(g))$. By symmetry the other inclusion holds.
The right-left implication can be shown with similar arguments.

Let us return to the proof of the lemma. The left-right implication
follows immediately from the claim and Lemma \ref{lemC:faceNew}.
For the right-left one, by using
the same arguments as in the proof of the previous lemma (i.e. the
boundedness of $\conv(E)$), we can show that $\supp(u)=\supp(u')$ and
$E_u(g)=E_{u'}(g)$, we then conclude by the claim above.
\end{proof}

Now we are ready to give the

\begin{proof}[Proof of Proposition \ref{propC:GrobnerCone_NC}]
By lemma \ref{lemC:NormalFace},
\[ \Norm_Q(\face_u(Q)) = \bigcap_{j=1}^r \Norm_{\New(g_j)}
(\face_u(\New(g_j))).\]
By the previous lemma,
\[\Norm_{\New(g_j)}(\face_u(\New(g_j)))= \{u' \in S(u) | \ini_u(g_j)=
\ini_{u'}(g_j)\}.\]
We then conclude by using Prop. \ref{propC:GrobnerCone_SB}.
\end{proof}

\subsection{Proof of Theorem \ref{theoC:SCGCpolyhedral}}

We recall that we are given an ideal $I$ in $\Ohat$. Let $u\in \Uloc$.

\begin{lem}\label{lemC:bah2.1}
For any $u'\in \overline{C_I[u]} \smallsetminus C_I[u]$ there exists a
local order $\prec$ such that the reduced standard bases of $I$ with
respect to $\prec_u$ and to $\prec_{u'}$ are the same.
\end{lem}

\begin{proof}
The proof is inspired by \cite[Prop. 2.1]{Compos}. We define the order
$\prec$ in a lexicographical way by $u$ and by any local order $<$. 
Note that $\prec=\prec_u$. Let $G$ be the $\prec$-reduced standard basis
of $I$. In order to prove the lemma, it is enough (see lemma
\ref{lemC:prec_12}) to show the following
for any $g\in G$: $\exp_{\prec_u}(g)= \exp_{\prec_{u'}}(g)$.
\begin{eqnarray*}
\exp_{\prec_{u'}}(g) & = & \exp_{\prec_{u'}}(\ini_{u'}(g))\\
 & = & \exp_{\prec_u}(\ini_{u'}(g)) \textrm{ by definition of } \prec.
\end{eqnarray*}
To finish it remains to prove this claim:
$\exp_{\prec_u}(\ini_{u'}(g))= \exp_{\prec_u}(g)$.\\
First since $\ND(\ini_{u'}(g)) \subset \ND(g)$, we have $\exp_{\prec_u}(g)
\succeq_u \exp_{\prec_u}(\ini_{u'}(g))$. Let us prove the reverse inequality.
First we prove:
\[(\star) \qquad \ord^{u'}(g)=u' \cdot \exp_{\prec_u}(g).\]
For any $\alpha \in \ND(g)$, and for any $u'' \in C_I[u]$,
$u'' \cdot \exp_{\prec_u}(g) \ge u'' \cdot \alpha$. Since $u'$ is
the limit of elements $u'' \in C_I[u]$, we obtain the same inequality for
$u'$ which proves $(\star)$. Therefore $\exp_{\prec_u}(g) \in
\ND(\ini_{u'}(g))$.
The desired inequality is then a direct consequence of the definition
of $\exp_{\prec_u}(\ini_{u'}(g))$. The claim and the lemma are proven.
\end{proof}

Finally, we can prove Theorem \ref{theoC:SCGCpolyhedral}.
\begin{proof}
The proof is exactly the same as that of \cite[Prop. 2.4]{sturmfels}.
For the sake of completeness let us give the main arguments.
For $u\in \Uloc$, let $u'$ be in $\overline{C_I[u]} \smallsetminus
C_I[u]$. Take an order $\prec$ as in the previous lemma and let $G$ be
the $\prec_u$-reduced standard basis of $I$. Then by Prop.
\ref{propC:GrobnerCone_NC}, $C_I[u]= \Norm_Q(\face_u(Q))$ and
$C_I[u']= \Norm_Q(\face_{u'}(Q))$ where $Q$ is the Minkowski sum of
$\New(g)$ for $g\in G$. By hypothesis on $u'$, $\face_u(Q)$ is a face
of the polyhedron $\face_{u'}(Q)$ thus $\overline{C_I[u']}$ is a face
of the closed convex cone $\overline{C_I[u]}$.

Now let us check the two axioms for being a fan. For axiom (1), let
$F$ be a face of the closure of some $C_I[u]$. Take any $u'$ in the
relative interior of $F$ then by the arguments above, $F$ is equal to
the face $\overline{C_I[u']}$.
For axiom (2), let $u,u'$ be any in $\Uloc$ and consider the closed
convex cone $P=\overline{C_I[u]} \cap \overline{C_I[u']}$. We have
seen that for any $u''\in P$, $\overline{C_I[u'']}$ is a face of both
$\overline{C_I[u]}$ and $\overline{C_I[u']}$. Thus $P$ is a (finite)
union of common faces. By convexity of $P$ this union is a singleton.
\end{proof}

\section{The closed Gr\"obner fan in $\hzoDhat$ or $\hzoDan$ is a
polyhedral fan}

\begin{theo}\label{theo:SCGCpolyhedral}
Let $I$ be a homogeneous ideal in $\hzoDhat$ or in $\hzoDan$ then its
closed Gr\"obner fan $\bar{\E}(I,\Uloc)$ is a rational polyhedral fan.
\end{theo}

The outline of the proof for this theorem is analogous to Theorem
\ref{theoC:SCGCpolyhedral}.
However, there are technical differences.
So, in this section we shall mainly prove the
results specific to the differential case.

\subsection{Division in the graded rings $\gr^w(\hzoDhat)$ and
$\gr^w(\hzoDan)$}

Let $\prec$ be a total monomial order on $\N^{2n}$ (or equivalently on
the terms $x^\alpha \xi^\beta$ where $\xi_i$ is a commutative variable
corresponding to $\pd_i$). Such an order shall be called local admissible
or simply \emph{admissible} if for any $i$, $x_i \prec 1$ and
$x_i \xi_i \succ 1$.

In the following we need a generalized version of the division theorem
of \cite{acg01} in $\hzoDhat$. Indeed we will prove it for any admissible
order. However we only treat the formal case and we will refer to
\cite{acg01} for the analytic case.

The next step will be a division theorem in the graded ring $\gr^w
(\hzoDhat)$ for any local weight vector $w\in \Wloc$.\\

Take an admissible order $\prec$ on $\N^{2n}$. Define the
order $\prec^h$ on $\N^{2n+1}$:

\[(\alpha, \beta, k) \prec^h (\alpha', \beta', k')
\iff
\begin{cases}
k+|\beta| < k' +|\beta'| \quad \textrm{ or} \\
k+|\beta| = k' +|\beta'| \textrm{ and }
 (\alpha ,\beta) \prec (\alpha', \beta').
\end{cases}
\]

As before we have the notion leading exponent $\exp_{\prec^h}$, leading
coefficient $\lc_{\prec^h}$, leading term $\lt_{\prec^h}$ and leading
monomial $\lm_{\prec^h}=\lc_{\prec^h} \cdot \lt_{\prec^h}$.

Let $P_1, \ldots, P_r \in \hzoDhat$. Consider the partition $\Delta_1 \cup
\cdots \cup \Delta_r \cup \bar{\Delta}$ of $\N^{2n+1}$ associated with the
$\exp_{\prec^h}(P_j)$ as done in the previous section.

\begin{theo}[Division theorem in $\hzoDhat$]
For any $P \in \hzoDhat$, there exists a unique
$(Q_1,\ldots,Q_r,R) \in (\hzoDhat)^{r+1}$ such that
\begin{itemize}
\item
$P=Q_1 P_1 + \cdots+ Q_r P_r +R$
\item
for any $j$, if $Q_j\ne 0$ then $\ND(Q_j) +\exp_{\prec^h}(P_i) \subset
\Delta_j$
\item
if $R\ne 0$, $\ND(R) \subset \bar{\Delta}$.
\end{itemize}
Moreover if $P$ and the $P_j$ are homogeneous then so are the $Q_j$ and $R$.
\end{theo}

In order to compare Gr\"obner fans in the analytic and the formal case, we
need a convergent division theorem in $\hzoDan$. This has been proved in
\cite{acg01}.

\begin{theo}[{\cite[Theorem 7]{acg01}}]\label{theo:acgTh7}
Let $w\in \Wloc$ and define the admissible order $\lhd$ in a lexicographical
way by $[|\beta|,\ <_0]$ where $<_0$ is any fixed admissible order then
define the admissible order $\lhd_w$ by refining $w$ by $\lhd$.
Then if we take $\prec=\lhd_w$ in the theorem above then the following holds:
if $P$ and the $P_j$ are in $\hzoDan$ then so are the $Q_j$ and $R$.
\end{theo}

Now let us prove the division theorem in $\hzoDhat$.
According to \cite{robbiano}, the order $\prec$ is defined by $2n$ independent
weight vectors. Let $w=(u,v)$ be the first of them. Since $\prec$ is
admissible, we have $u_i\le 0$ and $u_i+v_i \ge 0$, that is $w \in \Wloc$.
The first step of the proof of the division theorem is to reduce to the
case where $u$ has no zero components.

\begin{lem}
There exists a weight vector $w'=(u',v') \in \Wloc$ with $u'_i<0$ and
$u'_i+v'_i>0$ such that the order $\prec_{w'}^h$ defined first by
$k+|\beta|$, then by $w'$ and finally by $\prec$ satisfies:
$\exp_{\prec_{w'}^h}(P_j)= \exp_{\prec^h}(P_j)$ for any $j=1,\ldots,r$.
\end{lem}

For the proof of this lemma, we refer to \cite[Prop. 8]{acg01}.
Now once this lemma is proven then proving the division theorem for
$\prec_{w'}^h$ shall prove the theorem for $\prec^h$.

\begin{proof}[Proof of the formal division theorem]
The unicity is easy, and the statement concerning homogeneity also. Let us
prove the existence.
In order to simplify the notation, we will omit the subscript $\prec^h$.\\
Define the following sequences $P^{(i)}, Q_1^{(i)}, \ldots, Q_r^{(i)},
R^{(i)}$ as follows:
\begin{enumerate}
\item
Put $(P^{(0)}, Q_1^{(0)}, \ldots, Q_r^{(0)}, R^{(0)})=(P, 0, \ldots,0, 0)$.
\item
For $i\ge 0$, if $P^{(i)}=0$ then put
\[(P^{(i+1)}, Q_1^{(i+1)}, \ldots, Q_r^{(i+1)}, R^{(i+1)})=(P^{(i)},
Q_1^{(i)}, \ldots, Q_r^{(i)}, R^{(i)}).\]
\item
if $\exp(P^{(i)}) \in \bar{\Delta}$ then
\[(P^{(i+1)}, Q_1^{(i+1)}, \ldots, Q_r^{(i+1)}, R^{(i+1)})=
(P^{(i)}-\lm(P^{(i)}), Q_1^{(i)}, \ldots, Q_r^{(i)}, R^{(i)}+ \lm(P^{(i)})).\]
\item
if not, then $\exists ! j \in \{1,\ldots,r\}$
such that $\exp(P^{(i)}) \in \Delta_j$; put

$P^{(i+1)}=P^{(i)}- \frac{\lc(P^{(i)})}{\lc(P_{j})} \cdot
(x,\pd,h)^{\exp(P^{(i)})-\exp(P_{j})} \cdot P_{j}$,\\

$Q_{j}^{(i+1)}=Q_{j}^{(i)} + \frac{\lc(P^{(i)})}{\lc(P_{j})} \cdot
(x,\pd,h)^{\exp(P^{(i)})-\exp(P_{j})}$,\\

for $j'\ne j$, $Q_{j'}^{(i+1)}=Q_{j'}^{(i)}$,\\

$R^{(i+1)}=R^{(i)}$.
\end{enumerate}

By construction, $R^{(i)}$ tends to some $R \in \hzoDhat$ since
$R^{(i+1)}-R^{(i)}$ is one monomial for which the exponent does not belong
to $\ND(R^{(i)})$ and moreover the degree of $R^{(i)}$ is bounded (by that
of $P$). The same occurs for $Q_j^{(i)}$ that tends to some $Q_j \in
\hzoDhat$. To prove theorem, it remains to prove that $P^{(i)}$ tends
to $0$. For this purpose, we see $\hzoDhat$ as a free $\Ohat$-module and
we make it a topological space with the $(x)$-adic topology.
If we write $P^{(i)}=\sum_{\beta k} p^{(i)}_{\beta k}(x) \pd^\beta h^k$
then $P^{(i)}$ tends to $0$ if and only if each $p^{(i)}_{\beta k}(x)$
tends to $0$ in $\Ohat$. Now, by construction, we have $P=P^{(i)} +
\sum_j Q_j^{(i)} P_j + R^{(i)}$ and for $i>0$:

\begin{equation}\label{eq:P^i}
\exp(P) \succ \exp(P^{(i)}) \succ \exp(P^{(i+1)}).
\end{equation}

By this relation there exists $d\in \N$ such that all the $P^{(i)}$ have
degree $\le d$. So we can write $P^{(i)}= \sum_{(\beta, k)\in F}
p^{(i)}_{\beta k}(x) \pd^\beta h^k$ with $F$ finite.

For each $i$, all the $\exp_{\prec^h}(p^{(i)}_{\beta k}(x)\pd^\beta h^k)$
with $(\beta,k)\in F$ are pairwise distinct so $\exp_{\prec^h}(P^{(i)})$
equals the maximum of them, thus by relation (\ref{eq:P^i}) and by finiteness
of $F$, the following holds:

Fix any $(\alpha, k)\in F$ then for any $i_0$ there exists $i_1 \ge i_0$
such that: for any $i\ge i_1$, $\exp_\prec(p^{(i)}_{\beta k}(x) \pd^\beta
h^k) \prec \exp_\prec(p^{(i_0)}_{\beta k}(x) \pd^\beta h^k)$.\\

By the previous lemma, we may assume that $\prec$ is defined by weight
vectors where the first of them $w=(u,v)$ is such that $u_i<0$.\\
Now fix $(\alpha,k)\in F$. The preceding relation implies

$\ord^u(p^{(i)}_{\beta k}(x))\le \ord^u(p^{(i_0)}_{\beta k}(x))$.\\
But since $u$ has non zero components, the set of $\alpha$
such that $u \cdot \alpha$ equals some constant is a finite set. So in the
previous relation we cannot have an equality for all the $i\ge i_1$.
Thus their exists $i_2\ge i_1$ such that for any $i\ge i_2$, we have:

$\ord^u(p^{(i)}_{\beta k}(x))< \ord^u(p^{(i_0)}_{\beta k}(x))$.\\
With this final statement, it is easy to conclude that $p^{(i)}_{\beta k}(x)$
tends to $0$ in $\Ohat$ (again thanks to the fact that all the $u_i$ are
$<0$).
\end{proof}

Now we shall derive from the theorem above a \emph{division theorem in
the graded ring $\gr^w(\hzoDhat)$ or $\gr^w(\hzoDan)$}.\\

In order to better understand the division in these graded algebras,
let us describe the latter more precisely.

\begin{lem}\label{lem:gr}
For $w=(u,v)$, let us assume (for simplicity) that:
\begin{itemize}
\item
$u_i<0$ for $1\le i \le n_2$ with $u_i+v_i=0$ for $1\le i\le n_1$ and $u_i
+v_i>0$ for $n_1< i \le n_2$,
\item
$u_i=0$, $v_i>0$ for $n_2<i\le n_3$,
\item
$u_i=v_i=0$ for $n_3<i\le n$.
\end{itemize}
Then the graded ring $\gr^w(\hzoDhat)$ is canonically
isomorphic to
\[\k[[x_{n_2+1}, \ldots, x_n]][x_1,\ldots,x_{n_2}, \xi_{n_1+1},
\ldots, \xi_{n_3}, \pd_1, \ldots, \pd_{n_1}, \pd_{n_3+1},
\ldots, \pd_n][h]
\]
where $\xi_i$ is a commutative variable and
for $a(x)\in {\bf k}[[x_{n_2+1}, \ldots, x_n]][x_1,\ldots,x_{n_2}]$,
$[\pd_i, a(x)]= \frac{\pd a(x)}{\pd x_i} \cdot h$.
Concerning $\gr^w(\hzoDan)$, we have the same result replacing
$\k[[x_{n_2+1}, \ldots, x_n]]$ with $\C\{x_{n_2+1}, \ldots, x_n\}$.
\end{lem}

\begin{proof}
The proof consists of a simple verification.
\end{proof}

Let $w$ be in $\Wloc$. Given an admissible order $\prec$,
we define the (admissible) order $\prec_w$ by refining the partial
order defined by $w$.

The ring $\gr^w(\hzoDhat)$ is a bi-graded ring. There is a graduation by
the total degree in the $\pd_i$, $\xi_i$ and $h$. There is another graduation
associated with $w$. So in order to avoid confusions, a homogeneous element
for the second graduation shall be called $w$-homogeneous.

Let $P_1,\ldots,P_r$ be bihomogeneous in $\gr^w(\hzoDhat)$. Let
$\N^{2n+1}=(\cup_j \Delta_j) \cup \bar{\Delta}$ be the
partition associated with the $\exp_{\prec^h}(P_j)$.

\begin{cor}[Division theorem in $\gr^w(\hzoDhat)$]
For any bihomogeneous $P \in \gr^w(\hzoDhat)$, there exists a unique
$(Q_1,\ldots,Q_r,R) \in (\gr^w(\hzoDhat))^{r+1}$ made of bihomogeneous
elements such that
\begin{itemize}
\item
$P=Q_1 P_1 + \cdots+ Q_r P_r +R$
\item
for any $j$, if $Q_j\ne 0$ then $\ND(Q_j) +\exp_{\prec^h}(P_i) \subset
\Delta_j$
\item
if $R\ne 0$, $\ND(R) \subset \bar{\Delta}$.
\end{itemize}
\end{cor}

\begin{rem}
The division theorem in $\hzoDhat$ is a particular case of the
latter if we take $w=(0)$.
\end{rem}

\begin{proof}[Proof of the corollary]
There exist homogeneous elements
$P', P'_1, \ldots, P'_r$ in $\hzoDhat$ such that
$\ini_w(P')=P$ and $\ini_w(P'_j)=P_j$. Now let us consider the division
of $P'$ by the $P'_j$ w.r.t. the order $\prec_w^h$: $P'=\sum_j Q'_j P'_j
+R'$ with
\[(\star') \qquad \ND(Q'_j)+\exp_{\prec_w^h}(P'_j) \subset \Delta_j
\textrm{ and } \ND(R')\subset \bar{\Delta},\]
where the partition is associated with the $\exp_{\prec_w^h}(P'_j)$.
By homogeneity,
it follows that $\ord^w(P') \ge \ord^w(Q'_jP'_j)$ and $\ord^w(P')\ge
\ord^w(R')$. Set $Q_j=\ini_w(Q'_j)$ (resp. $R=\ini_w(R')$) if
the inequality above is an equality and $Q_j=0$ (resp. $R=0$)
otherwise. Then we have $P=\sum_j Q_j P_j +R$.

Now, by definition of $\prec_w^h$, for any homogeneous $Q\in \hzoDhat$,
$\exp_{\prec_w^h}(Q)=\exp_{\prec^h}(\ini_w(Q))$. This implies that
the partition of $\N^{2n+1}$ associated with the $\exp_{\prec_w^h}
(P'_j)$ is the same as that associated with the $\exp_{\prec^h}(P_j)$.

Moreover for any $Q\in \hzoDhat$, $\ND(\ini_w(Q)) \subset
\ND(Q)$ so the relation $(\star')$ becomes
\[\ND(Q_j)+\exp_{\prec^h}(P_j) \subset \Delta_j
\textrm{ and } \ND(R)\subset \bar{\Delta}.\]
We remark that the $Q_j$ and $R$ are
$w$-homogeneous and then bihomogeneous. The existence is proven.
The unicity can be proved easily.
\end{proof}

If we take $\prec=\lhd$ (see Th. \ref{theo:acgTh7}) then with the
same proof, we derive:

\begin{cor}[Division theorem in $\gr^w(\hzoDan)$ for $\lhd$]
For any bihomogeneous $P \in \gr^w(\hzoDan)$, there exists a unique
$(Q_1,\ldots,Q_r,R) \in (\gr^w(\hzoDan))^{r+1}$ made of bihomogeneous
elements such that
\begin{itemize}
\item
$P=Q_1 P_1 + \cdots+ Q_r P_r +R$
\item
for any $j$, if $Q_j\ne 0$ then $\ND(Q_j) +\exp_{\lhd^h}(P_i) \subset
\Delta_j$
\item
if $R\ne 0$, $\ND(R) \subset \bar{\Delta}$.
\end{itemize}
\end{cor}

\subsection{Standard bases in the graded ring $\gr^w(\hzoDhat)$}

The symbol $\prec$ still denotes any admissible order.
Let $J$ be a bihomogeneous ideal in $\gr^w(\hzoDhat)$ or in $\gr^w(\hzoDan)$
(and in the latter case, we take $\prec=\lhd$). We define
\[\Exp_{\prec^h}(J)=\{\exp_{\prec^h}(f) | 0\ne f\in J\}.\]
This set is stable by sums in ${\bf N}^{2n+1}$, thus by Dickson lemma:

\begin{defin}
There exists a finite set $G=\{g_1,\ldots,g_r\} \subset J$ such that
\[\Exp_{\prec^h}(J)= \bigcup_{j=1}^r (\exp_{\prec^h}(g_j)+
{\bf N}^{2n+1}).\]
Such a set is called a $\prec^h$-standard basis of $J$.
\end{defin}
\noindent
Given $G\subset J$, the following statements are equivalent:
\begin{itemize}
\item
$G$ is a bihomogeneous $\prec^h$-standard basis of $J$.
\item
For any $f\in \gr^w(\hzoDhat)$, $f\in J$ if and only if the remainder of
the division of $f$ by $G$ (w.r.t. $\prec^h$) is zero.
\end{itemize}

In the following lemma, we gather the results needed in the sequel. 
The generalizations of the division theorem in the previous subsection 
are used 
to prove the following analogous statements with Lemmas \ref{lemC:SB->in},
\ref{lemC:reducedSB}, \ref{lemC:prec_12}.

\begin{lem}
\begin{itemize}
\item
Given a bihomogeneous ideal $J$ in $\gr^w(\hzoDhat)$ (resp. in
$\gr^w(\hzoDan)$) and an admissible order $\prec$ (resp. $\prec=\lhd$),
a reduced standard basis w.r.t. $\prec^h$ exists and is unique.
Moreover it is bihomogeneous.
\item
Let $\prec_1$, $\prec_2$ be admissible orders. If $G$ is a
homogeneous (the reduced) $\prec_1^h$-standard basis of $J$ and
$\exp_{\prec_1^h}(g)= \exp_{\prec_2^h}(g)$ for any $g\in G$ then $G$
is a homogeneous (the reduced) $\prec_1^h$-standard basis of $J$.
\item
Given $I$ homogeneous in $\hzoDhat$ (resp. $\hzoDan$), if $G$ is the
reduced $\prec_w^h$-standard basis of $I$ then $\ini_w(G)$ is the
reduced $\prec^h$-standard basis of $\ini_w(I)$. Moreover we have
$\Exp_{\prec_w^h}(I)=\Exp_{\prec^h}(\ini_w(I))$.
\end{itemize}
\end{lem}


%

\subsection{Back to the Gr\"obner fan}

For $w=(u,v)\in \Wloc$, we call \emph{support of $w$} the set
$\supp(w) \subset \{1,\ldots,n\}^2$ defined as:
\[\underbrace{\{i| u_i<0\}}_{M(w)} \times
\underbrace{\{i|u_i+v_i>0\}}_{P(w)}. \]
Note that $M(w)$ and $P(w)$ are independent.
A consequence of lemma \ref{lem:gr} is:
For $w=(u,v) \in \Wloc$,
\[S(w)=S(w') \iff \supp(w)=\supp(w').\]

\begin{prop}\label{prop:GrobnerCone_SB}
Let $w\in \Wloc$. Given a homogeneous ideal $I$ in $\hzoDhat$ (resp. in
$\hzoDan$) and an arbitrary admissible order $\prec$ (resp. the order
$\prec=\lhd$), the $\prec^h_w$-reduced standard basis $G$ of $I$ satisfies:
\[C_I[w]=\{ w'\in \Wloc | \supp(w)=\supp(w') \textrm{ and }
\forall g\in G, \ini_w(g)=\ini_{w'}(g)\}.\]
\end{prop}

\begin{proof}
Exactly the same as that of Proposition \ref{propC:GrobnerCone_SB}.
\end{proof}

\begin{cor}\label{cor:FormAnalGF}
Let $I$ be a homogeneous ideal in $\hzoDan$ then the Gr\"obner fans
$\E(I,\Wloc)$ and $\E(\hzoDhat \cdot I, \Wloc)$ are equal.
\end{cor}

\begin{proof}
We put $\prec=\lhd$ and use the same arguments as for
Corollary \ref{corC:FormAnalGF}.
\end{proof}

After this corollary, we will work in $\hzoDhat$ for the rest of the
section.

Let us denote by $\pi$ the natural projection $\N^{2n} \times \N \to
\N^{2n}$, $(\alpha, \beta,k) \mapsto (\alpha, \beta)$.
For $g\in \hzoDhat$ let us define the Newton polyhedron of $g$
as the following convex hull:
\[\New(g)=\conv\Big(\pi(\ND(g)) + \{(\alpha, \beta)
\in \Z^{2n}| \forall (u,v) \in \Wloc,
(u,v)\cdot (\alpha,\beta) \le 0 \} \Big).\]

If we \emph{denote by $\Wloc^*$ the polar dual cone of $\Wloc$} then the
set under bracket is $\Wloc^* \cap \Z^{2n}$.
Let us characterize $\Wloc^*$: for $i=1,\ldots,n$, let $e_i \in \N^{2n}$
the vector having $1$ in its $i$th component and zero for the others;
and let $e'_i \in \Z^{2n}$, $e'_i=(0,\ldots, 0, -1, 0, \ldots,0, -1, 0
,\ldots, 0)$ with the $-1$ placed at position $i$ and $n+i$. Then
\[\Wloc^*=\big\{ \sum_{i=1}^n \lambda_i e_i +\sum_{i=1}^n \lambda'_i e'_i|
\ \lambda_i, \lambda'_i \ge 0\big\}.\]

Let $E(g) \subset \pi(\ND(g))$ be a finite subset such that:
$\pi(\ND(g))+ ((\oplus_1^n \N e_i) \times (\oplus_1^n \N e'_i))=
E(g)+((\oplus_1^n \N e_i) \times (\oplus_1^n \N e'_i))$.
This is possible by Dickson lemma. We have
\begin{equation}\label{eq:Newg}
\New(g)= \conv(E(g)) + \Wloc^*,
\end{equation}
which assures that $\New(g)$ is strictly speaking a polyhedron.\\

Let $w$ be in $\Wloc$, $\prec$ be any admissible order and $G$ be
the $\prec_w^h$-reduced standard basis of $I$. Let
\[Q=\sum_{g\in G} \New(g)\]
be the Minkowski sum of the Newton polyhedra of the $g$ in $G$.

\begin{prop}\label{prop:GrobnerCone_NC}
\[C_I[w]=\Norm_Q(\face_w(Q))\]
where $\Norm_Q(\face_w(\cdot))$ denotes the normal cone of
the face w.r.t. $w$.
\end{prop}

The proof will be based on the following lemmas.





\begin{lem}
For $g\in \hzoDhat$ and $w\in \Wloc$, let $E_w(g)$ be the set of
elements in $E(g)$ which are maximum for the scalar product with $w$ then
\begin{equation}\label{eq:faceNew}
\face_w(\New(g))=\conv(E_w(g)) + \big\{ \sum_{i\notin M(w)} \lambda_i e_i
+\sum_{i \notin P(w)} \lambda'_i e'_i| \ \lambda_i, \lambda'_i \ge
0\big\}.
\end{equation}
\end{lem}

\begin{proof}
It is easy to see that $\face_w(\Wloc^*)$ is the second member of the sum
in the RHS. Moreover for any polyhedra,
there is a general identity $\face_w(P+P')=\face_w(P)+ \face_w(P')$
so it suffices to show the equality $\face_w(\conv(E(g)))=
\conv(E_w(g))$. This equality follows easily from the fact that
the height of $\face_w(\conv(E(g)))$ and $\conv(E_w(g))$ w.r.t. $w$
is the same.
This fact can be shown as in the proof of lemma \ref{lemC:faceNew}
\end{proof}

Let $g_1,\ldots,g_r$ be in $\hzoDhat$. Put $Q_j=\New(g_j)$ and let
$Q=\sum_j Q_j$ be their Minkowski sum.

\begin{lem}\label{lem:NormalFace}
For any $w\in \Wloc$,
\[\Norm_Q(\face_w(Q)) = \bigcap_{j=1}^r \Norm_{Q_i}(\face_w(Q_j)).\]
\end{lem}

We can prove it by reducing to the case of polytopes of this lemma.

Here is the last lemma before the proof of proposition
\ref{prop:GrobnerCone_NC}.

\begin{lem}
Let $g\in \hzoDhat$ then for $w,w'\in \Wloc$:
\[\big[ \ini_w(g)=\ini_{w'}(g) \textrm{ and } S(w)=S(w') \big] \iff
\face_w(\New(g))= \face_{w'}(\New(g)).\]
\end{lem}

\begin{proof}
The proof is the same as that of lemma \ref{lemC:ini=face} where one has
just to replace $\N^n$ with $\Wloc^* \cap \Z^{2n}$.
\end{proof}

\begin{proof}[Proof of Proposition \ref{prop:GrobnerCone_NC}]
By lemma \ref{lem:NormalFace},
\[ \Norm_Q(\face_w(Q)) = \bigcap_{j=1}^r \Norm_{\New(g_j)}
(\face_w(\New(g_j))).\]
By the previous lemma,
\[\Norm_{\New(g_j)}(\face_w(\New(g_j)))= \{w' \in S(w) | \ini_w(g_j)=
\ini_{w'}(g_j)\}.\]
We then conclude by using Prop. \ref{prop:GrobnerCone_SB}.
\end{proof}

\subsection{Proof of Theorem \ref{theo:SCGCpolyhedral}}

Let $I$ be a homogeneous ideal in $\hzoDhat$.

\begin{lem}
Let $w' \in \overline{C_I[w]} \setminus C_I[w]$ for some $w\in \Wloc$. Then
there exists an admissible order $\prec$ such that the reduced standard
bases of $I$ w.r.t. $\prec_w^h$ and to $\prec_{w'}^h$ agree.
\end{lem}

\begin{proof}
The proof is a slight modification of \cite[Prop. 2.1]{Compos}. Since it
is similar to that of lemma \ref{lemC:bah2.1}, we sketch it.
We define the order $\prec$ first by $w$ and then by any admissible
order $<$. Let $G$ be the $\prec^h$-reduced standard basis of $I$.
Note that for homogeneous elements, $\exp_{\prec^h}$ coincides with
$\exp_{\prec_w^h}$.
Now, to prove the lemma, it is enough to show the following for any
$g\in G$: $\exp_{\prec_w^h}(g)=\exp_{\prec_{w'}^h}(g)$.
\begin{eqnarray*}
\exp_{\prec_{w'}^h}(g) & = & \exp_{\prec_{w'}^h}(\ini_{w'}(g))
 \textrm{ because } g \textrm{ is homogeneous}\\
& = & \exp_{\prec_w^h}(\ini_{w'}(g)) \textrm{ by definition of } \prec\\
& = & \exp_{\prec_w^h}(g).
\end{eqnarray*}
The last equality can be proved as in the proof of lemma \ref{lemC:bah2.1}.
\end{proof}

By \cite[Th. 4]{acg01}, we know that the number of different $\ini_w(I)$
is finite, so in order to finish the proof of Theorem
\ref{theo:SCGCpolyhedral} it remains to prove that $\bar{\E}(I,\Wloc)$
satisfies the two axioms for being a complex which can be done with
the same arguments as in the proof of Theorem \ref{theoC:SCGCpolyhedral}.

When considering applications or when doing calculations, we often need
to consider Gr\"obner fans restricted to linear subspaces in the space
of the weight vectors.

\begin{cor}
Let $L$ be a linear subspace in $\R^{2n}$ then for any homogeneous
ideal $I$ in $\hzoDhat$, the restriction of $\bar{\E}(I,\Wloc)$ to
$L$ is a polyhedral fan in $\Wloc \cap L$.
\end{cor}

\begin{proof}
Since we have proved that the local Gr\"obner fan is a polyhedral fan,
the restriction is also a polyhedral fan \cite[p. 195]{ziegler}.
\end{proof}

\section{Existence and comparison of the different Gr\"obner fans}

In what follows, we will discuss about the relation between the several
kinds of Gr\"obner fans. We will discuss both for polynomial rings and
associate rings and the case of differential operators.

Proofs are analogous each other. However statements are different
because the spaces of the weight vectors are different.

In order to clarify the difference, we present lemmas and propositions
for all the cases. The proofs will be made 
for the case of differential rings.\\

In the following $\Uloc'$ denotes the $u\in \Uloc$ such that $u_i<0$
and $\Wloc'$ denotes the $w=(u,v)\in \Wloc$ with $u_i<0$ and $u_i+v_i >0$.
In other terms, $\Uloc'$ and $\Wloc'$ are the interior of $\Uloc$ and
$\Wloc$ respectively.

\subsection{Without homogenization}

\subsubsection{The local Gr\"obner fan}

\

In order to compare the Gr\"obner fans, let us deal with the graded
algebras first.

\begin{rem}\label{rem:gr}
\
\begin{itemize}
\item
Take $u \in \Uloc$ and assume (for simplicity) that $u_i<0$ for
$1 \le i \le n_2$ and $u_i=0$ for $n_2 < i \le n$ then $\gr^u(\Ohat)$ is
canonically isomorphic to $\k[[x_{n_2+1}, \ldots, x_n]]
[x_1,\ldots, x_{n_2}]$.
\item
Take $w=\uv \in \Wloc$ and assume that the variables
$x_i$ and $\pd_i$ are ordered in a way that:
\begin{itemize}
\item
$u_i < 0$ for $1 \le i \le n_2$ with $u_i+v_i=0$ for $1 \le i \le n_1$
and $u_i+v_i >0$ for $n_1< i\le n_2$,
\item
$u_i=0$, $v_i >0$ for $n_2 < i \le n_3$,
\item
$u_i=v_i=0$ for $n_3 < i \le n$.
\end{itemize}
Then $\gr^w(\Dhat)$ is canonically isomorphic to
\begin{equation}\label{eq:identif}
\k[[x_{n_2+1}, \ldots, x_n]][x_1,\ldots, x_{n_2}, \xi_{n_1+1}, \ldots,
\xi_{n_3}, \pd_1,\ldots, \pd_{n_1}, \pd_{n_3+1},
\ldots, \pd_n]
\end{equation}
where $\xi_i$ is a commutative variable.
\end{itemize}
\end{rem}

In the following result, $\k[x]_{\langle x_{n_2+1},\ldots, x_n\rangle}$
shall denote the localization with respect to the prime ideal
generated by $x_{n_2+1},\ldots, x_n$ (or equivalently the localization
along the space defined by this ideal).
\begin{lem}\label{lem:grcompar}
Let $u$ be in $\Uloc$ and $w=\uv$ be in $\Wloc$. Let us take the same
situation as above, then:
\begin{itemize}
\item[(i)] $\dps \gr^u(\Ohat)= \k[[x_{n_2+1}, \ldots,
x_n]][x_1,\ldots,x_{n_2}] \bigotimes_{\k[x]} \gr^u(\Oalg)$.
\item[(ii)] $\dps \gr^u(\Oalg)=\k[x]_{\langle x_{n_2+1}, \ldots, x_n
\rangle} \bigotimes_{\k[x]} \gr^u(\k[x])$.
\item[(iii)] If $\k=\C$, we have two similar results, one for
$\gr^u(\Ohat)$ as a tensor of $\gr^u(\Oan)$ and the other one for
$\gr^u(\Oan)$ as a tensor of $\gr^u(\Oalg)$.
\item[(1)] $\dps \gr^w(\Dhat)=\k[[x_{n_2+1}, \ldots,
x_n]][x_1,\ldots,x_{n_2}] \bigotimes_{\k[x]} \gr^w(\Dalg).$
\item[(2)] $\dps \gr^w(\Dalg)=\k[x]_{\langle x_{n_2+1}, \ldots, x_n
\rangle} \bigotimes_{\k[x]} \gr^w(D)$.
\item[(3)] If $\k=\C$, we have two similar results, one for
$\gr^w(\Dhat)$ as a tensor of $\gr^w(\Dan)$ and the other one for
$\gr^w(\Dan)$ as a tensor of $\gr^w(\Dalg)$.
\end{itemize}
\end{lem}

Roughly speaking, this lemma says that the local variables $x_i$ for
which $u_i <0$ become global and the other ones stay unchanged. The
following is then trivial.

\begin{cor}
For any $u\in \Uloc'$, $\gr^u(\Ohat) = \gr^u(\Oan) = \gr^u(\Oalg) =
\gr^u(\k[x])$.
For any $w \in \Wloc'$, we have $\gr^w(\Dhat) = \gr^w(\Dan)= \gr^w(\Dalg)
= \gr^w(D)$.
\end{cor}

\begin{proof}[Proof of the Lemma]
For statement (1), the equality is trivial since $\gr^w(D)$ is
canonically isomorphic to $\k[x_1,\ldots, x_n, \xi_{n_1+1}, \ldots,
\xi_{n_3}, \pd_1,\ldots, \pd_{n_1}, \pd_{n_3+1},
\ldots, \pd_n]$.

Let us prove statement (2). First, let us remark that both rings in
this statement are subrings of $\gr^w(\Dhat)$ then it is enough to
prove the double inclusion. For the left-right one, it is enough to
show that $\ini_w(P)$ is in the RHS for any $P\in \Dalg$. Such a
$P$ is equal to $1/q(x) P'$ with $P'\in D$ and $q(x) \in \k[x]$ with
$q(0) \ne 0$. Since $\ini_w(1/q(x) P')=\ini_w(1/q(x))
\ini_w(P')$, let us show that $\ini_w(1/q(x))$ is in
$\k[x]_{\langle x_{n_2+1}, \ldots, x_n \rangle}$.

For simplicity, assume that $q(0)=1$ and write $q(x)=1-v_0-v_1$ where
$v_0=v_0(x_{n_2+1}, \ldots, x_n)$ and $v_1$ is in the ideal generated
by $x_1,\ldots,x_{n_2}$. As a series in $\k[[x]]$, we have
\[1/q(x)=V_1+ 1+v_0+v_0^2+ \cdots\]
where $V_1 \in \sum_{i=1}^{n_2} \k[[x]] x_i$. By
comparing the $w$-orders, we see that
$0=\ord_w(1+v_0+v_0^2+\cdots) > \ord_w(V_1)$. As a consequence
$\ini_w(1/q(x))=1+v_0+v_0^2+\cdots=1/(1-v_0)$ and it is in
$\k[x]_{\langle x_{n_2+1}, \ldots, x_n \rangle}$.

The right-left inclusion follows from the following: for any $P$ in
the RHS, $P$ is a finite sum of elements of the form $1/q \otimes
\ini_w(P')$ where $P'\in D$ and $q \in \k[x_{n_2+1}, \ldots, x_n]$
with $q(0)\ne 0$. But for such elements, we have $\ini_w(1/q)=1/q$
so $\ini_w(1/q P')=1/q \cdot \ini_w(P')$ and then $P$ belongs to
the LHS.

Finally, the proof of statement (3) uses the same arguments as that of
(1) and (2).
\end{proof}

Our main goal in this part is to present an algorithm for computing
two kinds of fan, the first one is the following.

\begin{defin}\label{def:locfan} \rm
Let $I$ be an ideal in $\k[x]$ (resp. in $D$). 
The set $\E(\Oalg I,\Uloc)$
(resp. $\E(\Dalg I, \Wloc)$) shall be called the local (open) Gr\"obner fan
of $I$ (at $0 \in \k^n$).
\end{defin}

Take $u\in \Uloc'$ (resp. $w \in \Wloc'$). We have seen that in this case,
all the graded rings that we considered are equal to $\gr^u(\k[x])$ (resp.
$\gr^w(D)$), so the initial ideals considered in the next proposition are
comparable.

\begin{prop}\label{prop:Wloc'}
Take $u\in \Uloc'$ and $w \in \Wloc'$.
\begin{itemize}
\item[(i)] If $I \subset \k[x]$ then $\ini_u(\Oalg I) = \ini_u(I)$.
\item[(ii)] If $I \subset \Oalg$ then $\ini_u(\Ohat I) = \ini_u(I)$.
\item[(iii)] If $\k=\C$ and $I\subset \Oalg$ then $\ini_u(\Oan I)=\ini_u(I)$.
\item[(iv)] If $\k=\C$ and $I\subset \Oan$ then $\ini_u(\Ohat I)=\ini_u(I)$.
\item[(1)] If $I \subset D$ then $\ini_w(\Dalg I) = \ini_w(I)$.
\item[(2)] If $I \subset \Dalg$ then $\ini_w(\Dhat I) = \ini_w(I)$.
\item[(3)] If $\k=\C$ and $I\subset \Dalg$ then $\ini_w(\Dan I) =
\ini_w(I)$.
\item[(4)] If $\k=\C$ and $I\subset \Dan$ then $\ini_w(\Dhat I) =
\ini_w(I)$.
\end{itemize}
\end{prop}

\begin{proof}
As we said, all these initial ideals are in a same ring $\gr^w(D)$.
Let us first suppose that $I$ is in $D$ and let us prove that
$\ini_w(\Dhat I) =\ini_w(I)$. Using the inclusions
\[\ini_w(I) \subset \ini_w(\Dalg I) (\subset \ini_w(\Dan I)) \subset \ini_w(\Dhat I),\]
claim (1) follows from the inclusion $\ini_w(\Dhat I) \subset
\ini_w(I)$, that we shall subsequently prove.

For this inclusion, it is enough to prove that for $P\in \Dhat I$,
$\ini_w(P) \in \ini_w(I)$. Such a $P$ can be written as a finite
sum
\[P=\sum_j c_j(x) P_j\]
with $c_j(x) \in \k[[x]]$ and $P_j \in I$. Put $m=\ord_w(P)$ and
for each $j$, let us write $c_j(x)=\sum_{\alpha} c_{j,\alpha}
x^\alpha$ where $c_{j,\alpha} \in \k$. Now let us decompose
$c_j(x)=p_j(x) + s_j(x)$ where
\[p_j(x) = \sum_{u \cdot \alpha + \ord_w(P_j) \ge m} c_{j,\alpha} x^\alpha.\]
Since all the $u_i$ are negative, each $p_j(x)$ is a polynomial. By
construction, we have $\ord_w(s_j(x) P_j) <m$. Therefore,
\[\ini_w(P)= \ini_w \left(\sum_j p_j(x) P_j \right) \in \ini_w(I),\]
and the inclusion is proved.

Finally, it is easy to remark that the method used for proving claim
(1) can be adapted with few changes to prove the other claims.
\end{proof}

As a trivial consequence, we obtain:

\begin{cor}\label{cor:Wloc'}
For $I$ in $\k[x]$, $\E(I,\Uloc')= \E(\Oalg I, \Uloc')$.
For $I \subset D$, $\E(I,\Wloc')= \E(\Dalg I, \Wloc')$.
\end{cor}

In other words, if we restrict ourselves to the interior of the
weight vector space then, global and local Gr\"obner fans do agree.
So the obstruction to have equality should be found in the border and
indeed they are not equal in general.

\begin{ex}\label{exC:border} \rm
\begin{enumerate}
\item
Take $I=\langle 1 \rangle$ in $\k[x]$. Then $\E(\Ohat I, \Wloc)$ is
equal to $\E(I,\Wloc) \cap S_\Ohat$ where $S_\Ohat$ is the stratification
given by the graded rings.
After this example we may ask if such a relation always happens.
We shall see in Th. \ref{theo:glob=locIFhomo} that this relation
is true if the ideal $I$ is homogeneous for a weight vector having
positive components. In the general case, the relation does not hold:
\item
Take $I=\langle 1-x_3, x_1+x_2 \rangle$ in $\k[x_1,x_2,x_3]$.
For the weights, we consider the space
$U_{12}=\{(u_1,u_2,0)|u_i\le 0 \}$.
In this space, the open global Gr\"obner fan has three cones $\{u_1 <u_2\}$,
$\{u_1=u_2\}$, $\{u_1>u_2\}$ while the local Gr\"obner fan is trivial
since $1-x_3$ is invertible. Here trivial means that we have four cones
$\{0\}, \{u_1<0, u_2=0\}, \{u_1<0, u_2<0\}, \{u_1=0, u_2<0\}$.
In this example, $\E(\Ohat I, U_{12})$ has 4 cones and $\E(I,U_{12}) \cap
S_\Ohat$ has 6 cones.
\end{enumerate}
\end{ex}

For a given algebraic ideal $I$, we called $\E(I,\Uloc)$ (resp. in $\Wloc$)
the local Gr\"obner fan of $I$, let us precise this term.

\begin{prop}\label{prop:locfan}
\
\begin{itemize}
\item[(i)] For $I \subset \Oalg$, $\E(I, \Uloc) = \E(\Ohat I, \Uloc)$.
\item[(ii)] If $\k=\C$ and $I \subset \Oan$ then $\E(I,\Uloc) =
 \E(\Ohat I, \Uloc)$.
\item[(1)] For $I \subset \Dalg$, $\E(I, \Wloc) = \E(\Dhat I, \Wloc)$.
\item[(2)] If $\k=\C$ and $I \subset \Dan$ then $\E(I,\Wloc) =
 \E(\Dhat I, \Wloc)$.
\end{itemize}
\end{prop}

This implies in particular that for a given $I$ in $D$, the following
Gr\"obner fans do agree: $\E(\Dalg I,\Wloc)$, $\E(\Dan I, \Wloc)$ (if
$\k=\C$), $\E(\Dhat I, \Wloc)$. Thus, if $\k=\C$ the local Gr\"obner
fan defined in \ref{def:locfan} agrees with the analytic standard fan
constructed by Assi et al.~\cite{acg01}.

Before giving the proof, let us state the last result of this
subsection.

\begin{theo}\label{th:compar_in}
Let $u$ be in $\Uloc$ and $w$ be in $\Wloc$.
\begin{itemize}
\item[(i)] For $I \subset \k[x]$, $\gr^u(\Oalg) \ini_u(I) = \ini_u(
\Oalg I)$.
\item[(ii)] For $I \subset \Oalg$, $\gr^u(\Ohat) \ini_u(I) =
\ini_u( \Ohat I)$.
\item[(iii)] When $\k=\C$, we have similar results concerning (a) $I
\subset \Oalg$ and $\ini_u( \Oan I)$ and (b) $I\subset \Oan$ and
$\ini_u(\Ohat I)$.
\item[(1)] For $I \subset D$, $\gr^w(\Dalg) \ini_w(I) = \ini_w(
\Dalg I)$.
\item[(2)] For $I \subset \Dalg$, $\gr^w(\Dhat) \ini_w(I) =
\ini_w( \Dhat I)$.
\item[(3)] When $\k=\C$, we have similar results concerning (a) $I
\subset \Dalg$ and $\ini_w( \Dan I)$ and (b) $I\subset \Dan$ and
$\ini_w(\Dhat I)$.
\end{itemize}
\end{theo}

\begin{rem}\label{rem:existLocGF}
By the previous theorem, for $I \subset \k[x]$, we have
\[\E(I,\Uloc) \textrm{ is a refinement of } \E(\Oalg I,\Uloc),\]
and for $I$ in $D$, we have:
\[\E(I,\Wloc) \textrm{ is a refinement of } \E(\Dalg I,\Wloc).\]
\end{rem}

\begin{proof}[Proof of Theorem \ref{th:compar_in}]
\begin{itemize}
\item[(1)] To prove the equality, it is enough to prove the right-left
inclusion, the opposite one being trivial. For this, we are reduced to
prove: for $P\in \Dalg I$, $\ini_w(P)$ is in the LHS. For such a $P$,
there exists $q(x) \in \k[x]$ with $q(0)\ne 0$ such that
$q(x) P \in I$, which implies that $\ini_w(q(x)) \ini_w(P) \in \ini_w(I)$.
Now $\ini_w(q(x))=\ini_u(q(x))$ and since the $u_i$ are non positive,
$\ini_u(q(x))$ is invertible in $\gr^w(\Dalg)$ (and not only in 
$\Dalg$).
\item[(2)] The proof for this statement is based on the homogenized
version of this theorem to be proved in the next subsection (see Th.
\ref{th:compar_inH}). As in (1), only the right-left inclusion is non
trivial. Let $P\in \Dhat \cdot I$. By lemma \ref{lem:hfactor} (proved
independently in the next subsection), $h_\zo(P) \in
\hzoDhat h_\zo(I)$. By theorem \ref{th:compar_inH}, $\ini_w(h_\zo(P))$
belongs to $\gr^w(\hzoDhat) \ini_w(h_\zo(I))$. We dehomogenize and
remark that $[\ini_w(h_\zo(P))]_{|h=1}= \ini_w(P)$ to conclude.
\item[(3)] (a) is a consequence of (b) and (2) while (b) is a direct
consequence of theorem \ref{th:compar_inH}.
\end{itemize}
\end{proof}

Using the previous theorem, let us give the

\begin{proof}[Proof of Prop. \ref{prop:locfan}]
Statement (1) follows from the following implication: for any $w$,
$w'$, if $\ini_w(\Dhat I) \subset \ini_{w'}(\Dhat I)$ then
$\ini_w(I) \subset \ini_{w'}(I)$. Let us prove it. By the
previous theorem, we are reduced to prove:
\[(\star) \qquad \gr^w(\Dhat) \ini_w(I) \subset \gr^w(\Dhat) \ini_{w'}(I) \Rightarrow \ini_w(I) \subset \ini_{w'}(I).\]
This is a direct consequence of the faithful flatness of $\k[[x]]$
over $\k[x]_{\langle x_1,\ldots,x_n \rangle}$ \cite[p. 62]{matsumura}.

\details{ Let us detail the argument. We take the notations of remark
\ref{rem:gr}. Let us denote by $\tilde{x}$ the set of variables
$x_{n_2+1}, \ldots, x_n$, and by $\k[\tilde{x}]_{\langle \tilde{x}
\rangle}$ the localization of $\k[\tilde{x}]$ w.r.t. the maximal ideal
generated by the set $\tilde{x}$. Now, let us first remark that the
LHS of $(\star)$ implies that for $w$ and $w'$, graded rings
are implicitly equal. The LHS of $(\star)$ takes place in
$\k[[\tilde{x}]] \otimes_{\k[\tilde{x}]} \gr^w(D)$ and the RHS in
$\k[\tilde{x}]_{\langle \tilde{x} \rangle} \otimes_{\k[\tilde{x}]}
\gr^w(D)$. Now let us work in the $\k[\tilde{x}]_{\langle \tilde{x}
\rangle}$-modules category. Then LHS of $(\star)$ is equivalent to
\[
\k[[\tilde{x}]] \bigotimes_{\k[\tilde{x}]_{\langle \tilde{x} \rangle}} \ini_w(I) \subset \k[[\tilde{x}]] \bigotimes_{\k[\tilde{x}]_{\langle \tilde{x} \rangle}} \ini_{w'}(I)
\]
which is equivalent to the exactness of
\[
0 \rightarrow \k[[\tilde{x}]] \bigotimes_{\k[\tilde{x}]_{\langle \tilde{x}
\rangle}} \ini_{w'}(I) \hookrightarrow
\Big( \k[[\tilde{x}]] \bigotimes_{\k[\tilde{x}]_{\langle \tilde{x} \rangle}} \ini_{w'}(I) \Big)
+
\Big( \k[[\tilde{x}]] \bigotimes_{\k[\tilde{x}]_{\langle \tilde{x} \rangle}} \ini_{w}(I)\Big) \rightarrow 0
\]
(where $\hookrightarrow$ means the inclusion) but
$(\k[[\tilde{x}]] \otimes \ini_{w'}(I))+(\k[[\tilde{x}]] \otimes \ini_{w}(I))
= \k[[\tilde{x}]] \otimes \big(\ini_{w'}(I) + \ini_{w}(I) \big)$
thus by faithful flatness of $\k[[\tilde{x}]]$ over $\k[\tilde{x}]_{\langle
\tilde{x} \rangle}$, the LHS of $(\star)$ is equivalent to the
exactness of
\[
0 \rightarrow \ini_{w'}(I) \hookrightarrow \ini_{w'}(I)
+ \ini_{w}(I) \rightarrow 0
\]
which means that $\ini_{w}(I) \subset \ini_{w'}(I)$, and we are
done.}

For statement (2), the proof is similar and based on the faithful
flatness of $\Ohat$ over $\Oan$ \cite[p. 62]{matsumura}.
\end{proof}

\subsubsection{On the passage from local to global}

We already know (see Remark \ref{rem:existLocGF}) that the global
Gr\"obner fan of $I$ is a refinement of the local one (at $0$) so we could
ask whether we have a passage from local to global.

Given an ideal $I$ in $D$, we considered the local Gr\"obner fan
$\E(\Dalg I,\Wloc)$ at $0\in \k^n$. Now, for any $x^0 \in \k^n$, we
can consider the local Gr\"obner fan of $I$ at the point $x^0$ which is
defined by considering $(u,v) \in \Wloc$ as weights on the variables
$x_i-x_i^0$ and $\pd_i$ (note that by the affine change of
coordinates $x \mapsto x'=x-x^0$, the derivations are not affected
i.e. $\frac{\pd}{\pd x'_i}= \frac{\pd}{\pd x_i}$). Another definition
is $\E(\Dalg I_{x \mapsto x-x^0}, \Wloc)$.

The following question is natural: is the global Gr\"obner fan $\E(I,
\Wloc)$ the common refinement of the local ones ? The next example
shows that the answer is no.

\begin{ex} \rm
Let $f=1+x_1+x_2 \in \k[x_1,x_2]$. To simplify, we will talk about
open Gr\"obner cones.
The restriction to $\Uloc$ of the global Gr\"obner fan is made of 4 cones:
$\{0\}$, $C_1=\{u_1 < 0, u_2=0\}$, $C_2=\{u_1<0, u_2<0\}$ and
$C_3=\{u_1=0, u_2 < 0\}$.
Now for $x_1^0,x_2^0 \in \k$, write
$f=1+x_1^0+x_2^0+(x_1-x_1^0)+(x_2-x_2^0)$.
In the case $1+x_1^0+x_2^0\ne 0$, $f$ is invertible in
$\k[x_1,x_2]_{\langle x_1-x_1^0,x_2-x_2^0\rangle}$, thus the local
Gr\"obner fan is trivial at such a point; that means that it is made
of the 4 cones $\{0\}, C_1,C_2,C_3$.
Now if $1+x_1^0+x_2^0= 0$ then the local Gr\"obner
fan at $(x_1^0,x_2^0)$ is made of 6 cones: $\{0\}$, $C_1$,
$\{u_1< 0, u_2< 0, u_1<u_2\}$, $\{u_1 < 0, u_2< 0, u_1=u_2\}$,
$\{u_1 < 0, u_2 < 0, u_1>u_2\}$ and $C_3$.

Thus the common refinement of the (two) local Gr\"obner fans is the
latter and it is different from the global one.
\end{ex}


\subsection{With homogenization}

Gr\"obner fans can be described explicitly by using the notion
of reduced Gr\"obner or standard bases. For an ideal $I$ in $\k[x]$
or in $D$ it is not possible, in general, to construct reduced Gr\"obner
bases when the order is not a global order, that is when we work with
weight vectors having some negative components. That explains why we need
to work with one more variable in homogenized rings. In these kind of
rings we can compute reduced Gr\"obner bases for any orders. Another
consequence of homogenization is that Gr\"obner fans of homogenized ideals
have convex cones (see \cite{acg00} and \cite{sst}).

In this short subsection, we state homogeneous analogues to the results
of the preceding section.
Most of the proofs are the same and will not be written.
First, here are two lemmas useful for the sequel.

\begin{lem}\label{lem:compute_hoo}
For a given ideal $I \subset D$, we can compute $h_\oo(I)$ as follows:
take any admissible (global) order $\prec$ and compute a
$\prec_{\oo}$-standard basis $G$ of $I$ then $h_\oo(I)$ will be
generated by the $h_\oo(P)$ for $P\in G$.
\end{lem}

\begin{proof}[Sketch of proof]
Let $P$ be in $I$. By definition of $G$, we can write $P=\sum_j Q_j
P_j$ with $P_j\in G$ and $Q_j \in D$ with $\ord_\oo(P)\ge \ord_\oo(Q_j
P_j)$. The homogenization then gives: $h_\oo(P)= \sum_j h^{\ord_\oo(P)
- \ord_\oo(Q_j P_j)} \cdot h_\oo(Q_j) \cdot h_\oo(P_j)$. Since the set
of the $h_\oo(P)$, $P$ running over $I$, generates $h_\oo(I)$, we are
done.
\end{proof}

\begin{lem}\label{lem:hfactor}
Given $I$ in $D$. We can compute generators of $h_\zo(I)$ in
$\hzoD$ by homogenizing a $\prec_\zo$-standard basis of $I$ and we have
\[\dps h_\zo(\Dhat I)= \hzoDhat h_\zo(I).\]
\end{lem}

\begin{proof}[Sketch of proof]
Let $G$ be a $\prec_\zo$-standard basis of $I$ then with the same
proof as above, we prove the first statement. Now by using Buchberger
$S$-criterion, $G$ is also a $\prec_\zo$-standard basis of $\Dhat I$
which implies that $h_\zo(G)$ also generates $h_\zo(\Dhat I)$.
\end{proof}

As it is easy to see, a homogeneous counterpart of Prop. \ref{prop:Wloc'}
holds from which follows a counterpart of Cor. \ref{cor:Wloc'} (the
proofs are the same and omitted).

\begin{cor}\label{cor:Wloc'Hom}
For a $\zo$-homogeneous ideal $I \subset \hzoD$,
$\E(I,\Wloc')= \E(\hzoDalg I, \Wloc')$.
\end{cor}

Here is the homogeneous version of Prop. \ref{prop:locfan}.
\begin{prop}\label{prop:locfanH}
\
\begin{itemize}
\item[(1)] For $I \subset \Dalg$, $\E(h_\zo(I), \Wloc) =
\E(h_\zo(\Dhat I), \Wloc)$.
\item[(2)] If $\k=\C$ and $I \subset \Dan$ then $\E(h_\zo(I),\Wloc) =
 \E(h_\zo(\Dhat I), \Wloc)$.
\end{itemize}
\end{prop}
Statement (2) is a direct consequence of \ref{cor:FormAnalGF} but it
can be proved as Prop. \ref{prop:locfan}.
The proof of (1) is almost the same as that of Prop. \ref{prop:locfan},
providing the fact that $h_\zo(\Dhat I)=\hzoDhat h_\zo(I)$ (this is
Lemma \ref{lem:hfactor}).

Finally, here is the homogeneous counterpart of Th. \ref{th:compar_in}.
\begin{theo}\label{th:compar_inH}
Let $w$ be in $\Wloc$.
\begin{itemize}
\item[(1)] For $I \subset D$, $\gr^w(\hzoDalg) \ini_w(h_\zo(I)) =
\ini_w( h_\zo(\Dalg I))$.
\item[(2)] For $I \subset \Dalg$, $\gr^w(\hzoDhat)
\ini_w(h_\zo(I)) = \ini_w( h_\zo(\Dhat I))$.
\item[(3)] When $\k=\C$, we have similar results concerning (a) $I
\subset \Dalg$ and $\ini_w( h_\zo(\Dan I))$ and (b) $I\subset \Dan$
and $\ini_w(h_\zo(\Dhat I))$.
\end{itemize}
\end{theo}

\begin{rem}\label{rem:compar_inH}
As a consequence of (1) of this theorem, we obtain: for $I\subset D$,
\[\E(h_\zo(I), \Wloc) \text{ refines } \E(h_\zo(\Dalg I),\Wloc).\]
\end{rem}

\begin{proof}[Proof of the Theorem]
Concerning (1), the proof is the same as for Th. \ref{th:compar_in}.
Let us prove (2). By \cite[Cor. 3.3]{got}
and by the previous lemma, $\Exp_{\prec_w^h}(h_\zo(I))= \Exp_{\prec_w^h}
(h_\zo(\Dhat I))$. Let $G$ be a homogeneous $\prec_w^h$-standard
basis of $h_\zo(I)$. The equality above implies that $G$ is also
a $\prec_w^h$-standard basis of $h_\zo(\Dhat I)$.
As a consequence, $\ini_w(G)$ generates $\ini_w(h_\zo(I))$
over $\gr^w(h_\zo(\Dalg))$ and also generates $\ini_w(h_\zo(\Dhat I))$
over $\gr^w(h_\zo(\Dhat))$. This implies the desired equality.

Statement (3)(a) is a consequence of (2) and (3)(b). For the latter,
consider a standard basis of $I$ w.r.t. the order $\lhd_w^h$ (see section
3) then the initial forms $\ini_w( \cdot )$ form a system of generators
of $\ini_w(h_\zo(I))$ and of $\ini_w(\hzoDhat h_\zo(I))$. Finally,
we remark that $\hzoDhat h_\zo(I) = h_\zo(\Dhat I)$ which can be
easily shown by homogenizing a $\zo$-standard basis of $I$.
\end{proof}

\section{Algorithms for local Gr\"obner fans}

In this section, we will focus on the following problems:\\
Given an ideal $I$ in $\k[x]$ (resp. in $D$), find an algorithm for
computing the following local Gr\"obner fans 
$\bar \E(\Ohat I, \Uloc \cap L)$
and 
$\bar \E(h_\zo(\Dhat I), \Wloc \cap L)$).
Here, $L$ is a linear subspace in a space of weights.
Our approach is based on the fact that the local Gr\"obner
fan can be refined by the global Gr\"obner fan of some homogeneous ideal,
the latter being computable.

Since the fans 
$\bar \E(\Ohat I, \Uloc )$ and
$\bar \E(h_\zo(\Dhat I), \Wloc )$
are polyhedral fans, it is enough for obtaining the Gr\"obner fan
restricted to the space $L$
to enumerate
the maximal dimensional cones in the whole weight space.
All lower dimensional Gr\"obner cones are obtained by taking faces
of the maximal dimensional cones.
We proved that local and global Gr\"obner cones agree
in maximal dimensional strata of the whole weight space.
Hence, if we ignore complexity of computation,
we do not need to consider the problem of enumerating
the local Gr\"obner cones in a weight space 
restricted to the linear subspace $L$. 
However, the cost of enumerating all the Gr\"obner cones
in the maximal dimensional strata in the whole weight space
is very high in general.
For example, for an ${\mathcal A}$-hypergeometric system 
associated to the matrix $A=(1,2,3)$, 
there are more than 1500 maximal dimensional cones.
Our implementation could not finish the enumeration in 3 days.
However, the small Gr\"obner fan, which is a fan restricted to
the linear subspace $u_i + v_i = 0$,
consists of only $7$ maximal dimensional cones.
This is the main reason why we restrict our weight space to a linear
subspace.
We note that local and global Gr\"obner cones do not agree in general
in the restricted weight space $\Wloc \cap L$.
See Example \ref{example:bs_example1} (see also Ex. \ref{exC:border}).

Finally, let us mention that Jensen has developed a
software package called Gfan \cite{gfan} (see also Fukuda et al.
\cite{fjt}), which can compute the (global) Gr\"obner fan of a polynomial
ideal. In \cite{fjt}, the authors propose a theory of Gr\"obner fans
for non-homogeneous ideals, but these fans are not ``local''. Indeed
it is easy to construct an example of two polynomial ideals
having the same Gr\"obner fan in the sense of \cite{fjt} and for which
the local Gr\"obner fans are different: in $\k[x_1,x_2]$, consider
$I_1=(g)$ and $I_2=(1+g)$ where $g=x_1+x_2+x_1 x_2^2 + x_1^2 x_2$.

\subsection{The commutative case}

Let $I$ be in $\k[x]$. We want an algorithm for $\E(\Ohat I, \Uloc)$.

\begin{theo}\label{theo:glob=locIFhomo}
Suppose there exists a weight vector $\alpha \in (\N_{>0})^n$ such that
$I$ is $\alpha$-homogeneous then
\[ \E(\Ohat I, \Uloc)= \E(I, \Uloc) \cap S_\Ohat.\]
\end{theo}
In other terms, local and global Gr\"obner fans coincide up to the
stratification by the ring $\Ohat$.
\begin{proof}
It is easy to see that $\bar{\E}(I,\Uglob)$ is a polyhedral fan.
Indeed the proof is the same as in the known case where $\alpha=(1,\ldots,
1)$ (see \cite{sturmfels}). It then implies
that $\bar{\E}(I, \Uloc) \cap \bar{\E}(\Ohat, \Uloc)$ is a fan.
Moreover, a polyhedral fan depends only on its maximal cones
and by Cor. \ref{cor:Wloc'}, local and global Gr\"obner fans coincide
on the interior of $\Uloc$ so $\bar{\E}(\Ohat I, \Uloc)=
\bar{\E}(I, \Uloc) \cap \bar{\E}(\Ohat, \Uloc)$. This last equality
is equivalent to the one we had to prove.
\end{proof}

Now suppose that $I$ is not homogeneous for any positive weight vector.
Choose a weight vector $\alpha \in (\N_{>0})^n$. Let $\{f_i, i\}$ be
a given set of generators of $I$.
Let $h$ be a new variable and let $I^{(h)}$ be the ideal of $\k[x,h]$
generated by the set of $h_\alpha(f_i)$, where $h_\alpha(f_i)$ is the
$\alpha$-homogenization of $f_i$.

\begin{prop}\label{prop:refine_poly}
$\E(I^{(h)}, \Uglob)$ refines $\E(I, \Uglob)$.
\end{prop}

\begin{proof}
Given $u$ and $u'$, we have to prove this implication:
\[\ini_u(I^{(h)})=\ini_{u'}(I^{(h)}) \Rightarrow \ini_u(I)=\ini_{u'}(I).\]
Suppose the LHS true. Let $f$ be in $I$ and let us prove that
$\ini_u(f) \in \ini_{u'}(I)$. Let us write $f= \sum_i q_i f_i$.
By homogenizing, we obtain: there exists $l$ and $l_i$ for any $i$
such that
\[h^l h_\alpha(f) = \sum_i h^{l_i} h_\alpha(q_i) h_\alpha(f_i)\]
so $h^l h_\alpha(f)$ belongs to $I^{(h)}$. By taking the initial form
w.r.t. $u$, we obtain by hypothesis: $\ini_u(h^l h_\alpha(f)) =
\sum_j r_j \ini_{u'}(g_j)$ with $g_j \in I^{(h)}$. It suffices to
set $h=1$ to obtain the desired relation. Thus this inclusion holds:
$\ini_u(I) \subset \ini_{u'}(I)$. The reverse inclusion holds by
symmetry which ends the proof.
\end{proof}

Combining this proposition with Rem. \ref{rem:existLocGF}, we obtain:
\begin{cor}
The fan $\E(I^{(h)}, \Uloc)$ refines $\E(\Oalg I, \Uloc)$.
\end{cor}

\begin{algo}[{\bf Computation of ${\bar \E}(\Oalg I, \Uloc \cap L)$}]
\label{algo:GFpoly} \ \rm \\
Input: an ideal $I$ in $\k[x]$. A linear subspace $L$ in $\R^n$. \\
Output: the local Gr\"obner fan ${\bar \E}(\Oalg I, \Uloc \cap L)$. \\
Step 1:
\begin{quote}
Compute the set $\Sigma_0$ of maximal cones of the global Gr\"obner fan
of $I^{(h)}$ restricted to $\Uloc \cap L$: $\E(I^{(h)}, \Uloc \cap L)$.
\end{quote}
Step 2:
\begin{quote}
For any $C, C' \in \Sigma_0$, using an \'ecart division, 
compare $\ini_u(\Oalg I)$ and $\ini_{u'}(\Oalg I)$ 
in $\gr_u(\Oalg)$ 
for some $u\in C$, $u' \in C'$. 
If it is equal, glue $C$ and $C'$. 
By continuing this process, 
we construct the set $\Sigma$ of maximal cones of 
${\bar \E}(\Oalg I,\Uloc \cap L)$. 
From $\Sigma$, construct the set $\E$ of all the cones.
\end{quote}
Output $\E$.
\end{algo}

Let us add two remarks to the algorithm.

The step 1 can be performed by flipping of the maximal dimensional
cones in the linear space $L$ with respect to the facets.
As to details, see Algorithm 3.6 of the book by Sturmfels \cite{sturmfels}.
Note that the correctness of the flipping procedure for the enumeration
comes from the fact that the Gr\"obner fan is a polyhedral fan.
The flipping procedure can be accelarated by Collart-Kalkbrener-Mall's
Gr\"obner walk method \cite{ckm}.
We note that the method also works in rings of differential operators.

The second remark is on the \'ecart division in Step 2.
Suppose that $\supp\,u=\supp\,u' = \{1, 2, \ldots, m\}$.
Then, we have
$$ \gr^u (\Oalg) \simeq \left\{
  f/g  |  f \in \k[x], g \in \k[x_{m+1}, \ldots, x_n], g(0) \not= 0 
\right\}.
$$
The \'ecart division in literatures usually suppose that 
denominators are in $\k[x]$,
but our denominator $g$ lies in $\k[x_{m+1}, \ldots, x_n]$.
Hence, in our case, a reducer must be chosen
so that the multiplier monomial for the pseudo-division
is in the ring $\k[x_{m+1}, \ldots, x_n]$ 
instead of $\k[x]$ as in the case of the usual \'ecart division
\cite{mora}, \cite{grabe}, \cite{gp}.

\subsection{The non homogeneous case}

\begin{prop}
Let $I$ be in $D$ then $\E(h_\oo(I), W)$ refines $\E(I,W)$ for any
$W\subset \Wglob$.
\end{prop}

In fact, we have a more general result. Let $G$ be any system of
generators of $I$ over $D$ and let $I^{(h)} \subset h_\oo(D)$ be
generated by $\{h_\oo(P) | P\in G\}$ (of course, this ideal is not
uniquely determined).

\begin{rem}[On the computation of $\E(I^{(h)}, \Wloc)$] \rm
Since $I^{(h)}$ is homogeneous, it is well known that its Gr\"obner
fan can be computed by using reduced Gr\"obner bases w.r.t. to
well-orders (see e.g. \cite{sst}).
\end{rem}

\begin{prop}
The fan $\E(I^{(h)}, W)$ refines $\E(I,W)$ (for $W\subset \Wglob$).
\end{prop}

\begin{proof}
The same as that of Prop. \ref{prop:refine_poly}.
%
\end{proof}

Now combining the previous proposition with Remark
\ref{rem:existLocGF}, we obtain:

\begin{cor}\label{cor:refine_loc}
The fan $\E(I^{(h)}, \Wloc)$ is a refinement of $\E(\Dalg I, \Wloc)$.
\end{cor}

By this Corollary, we obtain the following algorithm.

\begin{algo}[{\bf Computation of $\E(\Dalg I, \Wloc)$}]\label{algo:GF}
\ \\
\rm
Input: an ideal $I$ in $D$.\\ 
Output: the local Gr\"obner fan $\E(\Dalg I, \Wloc)$. \\
Step 1:
\begin{quote} 
Compute one of the two (global) Gr\"obner fans and call it
$\E_0$:\\ 
$\bullet$ $\E(h_\oo(I), \Wloc)$ (see Lemma \ref{lem:compute_hoo} 
for how to compute $h_\oo(I)$)\\ 
$\bullet$ $\E(I^{(h)}, \Wloc)$\\ 
This computation can be done as in \cite[Chapter 2]{sst}.
\end{quote}
Step 2:
\begin{quote} 
Two cones $C$ and $C'$ of $\E_0$ are
said to be in a same class if for one $w \in C$ and one
$w' \in C'$, we have $\ini_w(\Dalg I)=\ini_{w'}(\Dalg
I)$.\\ 
By using Algorithm \ref{algo:compare_in}, compute the classes
of $\E_0$.\\ 
Set $\E:=$the set of the classes of $\E_0$.
\end{quote}
Output $\E$.
\end{algo}

We note that 
${\bar \E}(\Dalg I, \Wloc)$ is not a polyhedral fan in general.
Hence, we need to perform the merging procedure in Step 2 
for all dimensional Gr\"obner cones.

\begin{algo}\label{algo:compare_in} \ \\
\rm
Input: $I\subset D$, $w, w' \in \Wloc$\\ 
Output : $1$ if $\ini_w(\Dalg I)=\ini_{w'}(\Dalg I)$ and $0$ if not.\\ 
(1) Compute $G_1$ a $w$-standard basis of $I$ 
and $G_2$ a $w'$-standard basis of $I$.\\ 
(2) By a reduction via an \'ecart
(division as in \cite{got}), compare $G_1$ and $G_2$
\end{algo}

Let us make some remarks on Algorithm \ref{algo:compare_in}:
Concerning step (1), the following fact is basic: Let $H_1$
(resp. $H_2$) be the reduced Gr\"obner basis of $I^{(h)}$ w.r.t. a
well order that privileges $w$ (resp. $w'$). Then by the
specialization $h=1$ we can set $G_i={H_i}_{|h=1}$. Thus Step (1) does
not require extra computations since the reduced Gr\"obner basis of
$I^{(h)}$ were needed in Algorithm \ref{algo:GF}.

For Step (2), we can use two methods. The first one is an \'ecart
division in $D$ as in \cite{got} (which works although the \'ecart
division of loc. cit. is stated in $\hzoD$). The second method is
based on the following equivalence
\[\ini_w(\Dalg I)=\ini_{w'}(\Dalg I) 
 \iff (\ini_w(\Dalg I))^{(h)}= (\ini_{w'}(\Dalg I))^{(h)}.\]
The right-left implication can be proved in a similar way as that of
Prop. \ref{prop:refine_poly}, while the left-right one is
trivial. Using this equivalence, we can use \'ecart division in
$\hzoD$. 

\subsection{The doubly homogenized Weyl algebra}

In the next subsection, we present two variants for an algorithm of
computing the Gr\"obner fan of $h_\zo(\Dhat I)$.
One variant is based on the doubly homogenized Weyl algebra that
we introduce here.

The doubly homogenized Weyl algebra  $h'(D)$
is generated by
$$ x_1, \ldots, x_n ,  \pd_1, \ldots, \pd_n, h, h' $$
with the relations
\[  \pd_i x_i = x_i \pd_i + h h' .\]

Let $\prec$ be a total order on the set of normally ordered
monomials $x^\alpha \pd^\beta h^p {h'}^q$ in $h'(D)$.
Such an order is called a 
\emph{multiplicative monomial order\/} if the following two conditions
hold: \\
1. $x_i \pd_i \succ h h' $ for $i = 1, 2, \ldots, n$; \\
2. $x^\alpha \pd^\beta h^k {h'}^l \prec  
    x^{\alpha'} \pd^{\beta'} h^{k'} {h'}^{l'}$ 
   $\Longrightarrow$
   $x^{\alpha+a} \pd^{\beta+b} h^{k+q} {h'}^{l+r} \prec  
    x^{\alpha'+a} \pd^{\beta'+b} h^{k'+q} {h'}^{l'+r}$.
Under this definition, the theory of Gr\"obner basis 
works analogously with $D$ for this ring. 
The Gr\"obner fan is also defined analogously.
We denote by $t$ a weight for $h$ and by $t'$ a weight for $h'$.
\begin{theo}
Consider the weight space
$$ u_i + v_i \geq t + t' $$
and fix a homogeneous left ideal $J$ in $h'(D)$.
The collection of closures of the Gr\"obner cones of $J$
in the weight space is a polyhedral fan.
\end{theo}

The proof is analogous to the case of $\hooD$ (see \cite{sst}). Indeed
with the definitions below we dispose of a well order and
reduced Gr\"obner bases exist. As for Prop. \ref{prop:refine_poly},
we have:

\begin{prop}\label{prop:refine_doublyWeyl}
If we restrict the fan to the linear subspace $t=t'=0$,
then it is a refinement of the Gr\"obner fan of $J_{|h'=1}$ in $\hzoD$
\end{prop}

Let $\prec_1$ be a multiplicative monomial order in 
$\hzoD$ satisfying {\rm (2.1), (2.2), (2.3)} of \cite{got}:
\begin{enumerate}
\item  $x_i \pd_i \succ_1 h$.
\item if $|\beta|+k < |\beta'|+k'$ then $x^\alpha \pd^\beta h^k
\prec_1 x^{\alpha'} \pd^{\beta'} h^{k'}$.
\item  $x^\alpha \preceq 1$.
\end{enumerate}
We define an order on $h'(D)$ as a block order on $\prec_1$ as follows:
%
\begin{equation}\label{eq:dorder}
x^\alpha \pd^\beta h^k {h'}^l \prec x^{\alpha'} \pd^{\beta'} h^{k'} {h'}^{l'}
\Leftrightarrow
\begin{cases}
|\alpha|+|\beta|+k+l < |\alpha'|+|\beta'|+k'+l' \\
\textrm{or } =   \textrm{ and }
 x^\alpha \pd^\beta h^k \prec_1 x^{\alpha'} \pd^{\beta'} h^{k'}.
\end{cases}
\end{equation}

\begin{lem} \label{lemma:dehomo}
Let $f, g \in {h'}(D)$ be homogeneous
and $\prec$ be an order satisfying {\rm (\ref{eq:dorder})}.
If $\exp_{\prec}(f) \preceq \exp_{\prec}(g)$,
then $\exp_{\prec_1}(f_{|h'=1}) \preceq_1 \exp_{\prec_1}(g_{|h'=1})$.
\end{lem}
This lemma does not hold for $h' \mapsto h$.

\begin{proof}
Suppose that $\lm_\prec(f)$ is equal to
$ c x^\alpha \pd^\beta h^k {h'}^l$.
Since $f$ is homogeneous and $\prec$ is a block order, 
the leading monomial of $f_{|h'=1}$ is 
$\lm_{\prec_1}(f_{|h'=1}) = c x^\alpha \pd^\beta h^k$.
In other words, cancelation does not happen when $h'\mapsto 1$.
It completes the proof.
\end{proof}

\begin{theo}
Fix a homogeneous left ideal $J$ in $h'(D)$ generated by $F$.
Let $\prec_1$ be an ordering in $h_\zo(D)$
satisfying {\rm (2.1), (2.2), (2.3)} of \cite{got}:
Let $G$ be a Gr\"obner basis of $I$ for the order $\prec$
defined as in {\rm (\ref{eq:dorder})}.
Then, the dehomogenization $G_{|h'=1}$
is a standard basis in the sense of \cite[Def. 3.1]{got}
of the ideal generated by $F_{|h'=1}$ with respect to the order $\prec_1$.
\end{theo}

\begin{proof}
$G_{|h'=1}$ generates the ideal generated by $F_{|h'=1}$
since $G$ generates the ideal generated by $F$.
Denote the elements of $G$ by $\{ g_i \}$.
Since $G$ is a Gr\"obner basis, we have Buchberger's S-pair criterion;
\[ {\rm sp}(g_i, g_j) = \sum_k q_{ijk} g_k, \quad
   \exp_{\prec}({\rm sp}(g_i, g_j)) \succeq \exp_{\prec}(q_{ijk} g_k).\]
It follows from Lemma \ref{lemma:dehomo} that
${\rm sp}(g_i, g_j)_{|h'=1} = \sum_k {q_{ijk}}_{|h'=1} {g_k}_{|h'=1}$ with
\[\exp_{\prec_1}({\rm sp}(g_i, g_j)_{|h'=1}) \succeq_1
     \exp_{\prec_1}((q_{ijk} g_k)_{|h'=1}).\]
Then, by Theorem 3.2 of \cite{got}, we conclude the theorem.
\end{proof}

\subsection{The homogeneous case}

Denote by $J_1$ the ideal of $h'(D)$ generated by 
homogenization of $h_\zo(I)$ with the variable $h'$
(in fact we homogenize any set of generators of $h_\zo(I)$).
The ideal $J_1$ is homogeneous 
in the doubly homogenized Weyl algebra $h'(D)$. 
Combining Proposition \ref{prop:refine_doublyWeyl} 
and Remark \ref{rem:compar_inH}, we have:

\begin{cor}
The fan $\E(J_1,\Wloc)$ is a refinement of $\E(h_\zo(\Dalg I), \Wloc)$.
\end{cor}

This Corollary gives a first variant for computing local Gr\"obner fan.
Let us explain the second variant.
Now, we use the notions of \cite{got}: for $P\in \hzoD$, we denote by
$P^{(s)}$ the $(\mathbf{-1,1})$-homogenization in $\hzoD[s]$ (here $s$
is a new variable commuting with the other ones).
Let $G$ be a set of generators of $h_\zo(I)$ over $\hzoD$, then
define $h_\zo(I)^{(s)}$ as the ideal of $\hzoD[s]$ generated by
$P^{(s)}$ for $P\in G$. To simplify denote by $J_2$ this ideal.

Then we can consider the Gr\"obner fan $\E(J_2,\Wloc)$ as we did before,
by putting a weight $0$ on $s$.

\begin{rem}[On the computation of $\E(J,\Wloc)$]\label{rem:computeGFJ} \rm
In $\hzoD[s]$ we dispose of a well-order $\prec_s$ (notation of
\cite{got}), thus the Gr\"obner fan of $J$ can be computed via reduced
Gr\"obner bases computations.
\end{rem}

\begin{prop}\label{prop:refine_locH}
The fan $\E(J_2,\Wloc)$ is a refinement of $\E(h_\zo(I), \Wloc)$.
\end{prop}
We omit the proof since it is similar to that of
Prop. \ref{prop:refine_poly}. By Prop. \ref{prop:refine_locH} and
Remark \ref{rem:compar_inH}:

\begin{cor}\label{cor:refine_locH}
The fan $\E(J_2,\Wloc)$ is a refinement of $\E(h_\zo(\Dalg I), \Wloc)$.
\end{cor}

As a consequence, we obtain the next algorithm.

\begin{algo}[{\bf Computation of ${\bar \E}(h_\zo(\Dalg I), \Wloc \cap L)$}]\label{algo:GFH} \rm
\ \\
Input: an ideal $I$ in $D$. A linear subspace $L$ in $\R^{2n}$. \\
Output: the local Gr\"obner fan ${\bar \E}(h_\zo(\Dalg I), \Wloc \cap L)$. \\
Step 1:
\begin{quote}
Compute the set $\Sigma_0$ of the maximal cones 
of the global Gr\"obner fan ${\bar \E}(J, \Wloc \cap L)$ 
where $J$ is one of $J_1,J_2$.
\end{quote}
Step 2:
\begin{quote}
As before, compute the classes of $\Sigma_0$, 
by using Algorithm \ref{algo:compare_inH}.\\
From the set of the classes of $\Sigma_0$,
construct the set $\E$ of all the cones of $\E(J, \Wloc \cap L)$.
\end{quote}
Output $\E$.
\end{algo}

\bigbreak

\begin{algo}\label{algo:compare_inH} \ \rm \\
Input: $I\subset D$, $w, w' \in \Wloc$\\ 
Output : $1$ if $\ini_w(h_\zo(\Dalg I))=\ini_{w'}(h_\zo(\Dalg I))$ 
and $0$ if not.\\
(1) Compute $G_1$ a $w$-standard basis of $I$ and $G_2$ a
$w'$-standard basis of $h_\zo(I)$.\\
(2) By a reduction via an \'ecart division in $\gr^w(\hzoD)$ as
in \cite{got}, compare $G_1$ and $G_2$.
\end{algo}

Let us give some remarks on this algorithm.\\
- As before, we can obtain $G_1$ and $G_2$ by dehomogenizing operators
computed in Algorithm \ref{algo:GFH}. The justification of Step (1) lies
in the fact that $h_\zo(I)$ generates $h_\zo(\Dalg I)$ over $\hzoDalg$.
Moreover we have seen in Lemma \ref{lem:hfactor} how to compute $h_\zo(I)$.\\
- The \'ecart division of \cite{got} (stated in $\hzoD$) 
works well in $\gr^w(\hzoD)$ with a modification;
a reducer must be chosen so that the multiplier monomial 
for the pseudo division is in the ring $\k[x_{m+1}, \ldots, x_n]$ 
instead of $\k[x]$ as in the case of the usual \'ecart division.  
Here, we assume that ${\rm Supp}(u)=\{1, \ldots, m\}$. \\
- Finally, since $\bar{\E}(h_\zo(\Dhat I)),\Wloc$ is a polyhedral fan,
it is enough, in step 1, to compute only the maximal cones.

\subsection{Tips for implementation and examples}

We implemented the algorithms except Algorithm \ref{algo:GF} in this paper 
in a combination of Kan \cite{kan} and Polymake \cite{polymake}.
The following is a list of known implementations 
for Gr\"obner fans
and some efficiency studies including
Gr\"obner walk.
\begin{enumerate}
\item 
Macaulay \cite{mcly} command {\tt hull} implemented by A. Reebs
constructs state polytopes.
\item The Singular developing team \cite{singular} implements 
the Gr\"obner walk techniques \cite{ckm}, \cite{fglm}.
\item B. Huber and R. Thomas gave an efficient algorithm and implementation
specialized to getting Gr\"obner fan for toric ideals \cite{ht}.
\end{enumerate}
Our implementation is new in view of the following aspects.
\begin{enumerate}
\item It is the first implementation for local Gr\"obner fans.
\item Computation is split into a system for algebra (Kan)
and a system for geometry (Polymake).
They are connected by the OpenXM-RFC 104 protocol 
(OoHG, OpenXM on HTTP/GET) \cite{openxm}.
This design gives us a robust system, since the polymake is a strong,
flexible, and robust system for polytopes developed by E. Gawrilow and
M. Joswig linked with GMP and cdd by K. Fukuda.
We use polymake properties {\tt DIM}, {\tt FACETS},
{\tt INEQUALITIES}.
\end{enumerate}

We have explained our algorithms for local Gr\"obner fans
in the previous sections.
The dominant part in the computation time is Step 1's.
Let us explain some details of Step 1, which obtains the maximal dimensional
Gr\"obner cones in a homogenized ring.
\begin{enumerate}
\item Parametrize the weight space $\Wloc \cap L$ by a polyhedron
$W \subseteq {\R}^{{\rm dim}\, (\Wloc \cap L)}$.
The parameterization is given by a matrix $W_{\rm cone}$.
\item Find a starting weight vector $u$ such that the dimension
of the (projected) Gr\"obner cone $C[u] \subseteq W$ is
${\rm dim}\, (\Wloc \cap L)$.
\item  If $C[u]$ is not pointed, then there exists a non-zero linear space
in $C[u]$.
Let $L_{\rm cone}$  be the maximal linear space in $C[u]$.
We construct Gr\"obner cones in $W/L_{\rm cone}$. In other words
we work in the othogonal space to $L_{\rm cone}$.
It is necessary because the Polymake and polyhedral algorithms works 
efficiently only for pointed cones.
Note that all maximal dimensional cones in the $W$ space contain
the linear space $L_{\rm cone}$ since the Gr\"obner fan is a polyhedral fan.
\item Enumerate the facets of $C[u]/L_{\rm cone}$.
They are either on the border of $W/L_{\rm cone}$ or not.
\item Perturb the weight vector $u$ with respect to a facet $F$, 
which is not on the border of the weight space $W$,
as explained in \cite[Chapter 3]{sturmfels}.
\item Lift the new weight vector in $W/L_{\rm cone}$ to $\Wloc \cap L$.
Let $u'$ be the lifted weight vector.
Construct the new reduced Gr\"obner basis and construct $C[u']/L_{\rm cone}$.
If $C[u]/L_{\rm cone} \cap C[u']/L_{\rm cone}$  is the facet $F$,
then we mark the facet $F$ as flipped,
else retry to get a new $u'$ with a smaller $\varepsilon$ perturbation.
\item Continue the procedure of getting new Gr\"obner cones and new facets
until all facets have been marked with flipped.
\end{enumerate}
As to more details, see the source code {\tt gfan.sm1} 
of our test implementation.

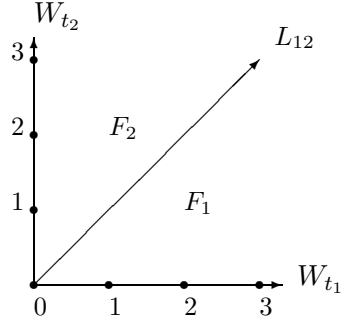
\begin{figure}[htb]  \label{figure:bs_example1} 
\setlength{\unitlength}{1mm}
\begin{picture}(45,45)
\put(10,10){\circle*{1}} \put(10,6){$0$}  
\put(20,10){\circle*{1}} \put(20,6){$1$}
\put(30,10){\circle*{1}} \put(30,6){$2$}
\put(40,10){\circle*{1}} \put(40,6){$3$}
\put(10,20){\circle*{1}} \put(7,20){$1$}
\put(10,30){\circle*{1}} \put(7,30){$2$}
\put(10,40){\circle*{1}} \put(7,40){$3$}
\put(10,10){\vector(1,0){33}}
\put(45,10){{\large $W_{t_1}$}}
\put(10,10){\vector(0,1){33}}
\put(10,45){{\large $W_{t_2}$}}
\put(10,10){\vector(1,1){30}}
\put(42,42){$L_{12}$}

\put(30,20){$F_1$}
\put(20,30){$F_2$}
\end{picture}
\caption{Gr\"obner fan for Example \ref{example:bs_example1}}
\end{figure}

\begin{ex} \label{example:bs_example1} \rm 
Let us consider the left ideal $I$ in 
$D={\bf Q} \langle t_1, t_2, x, y, \pd_{t_1}, \pd_{t_2}, \pd_x, \pd_y \rangle$
generated by 
\begin{equation}  \label{eq:bs_example1}
g_1=t_1-y,
g_2=t_2-(y-(x-1)^2),
g_3=(-2x+2)\partial_{t_2}+\partial_{x},
g_4=\partial_{t_1}+\partial_{t_2}+\partial_{y}.
\end{equation}
The ideal is used for computing a Bernstein-Sato polynomial
for $f_1=y$ and $f_2 = y-(x-1)^2$ via a local Gr\"obner fan
\cite{Compos}.

We will examine several Gr\"obner fans in the weight space
$$ W = \{ (-w_1,-w_2, 0,0, w_1, w_2, 0,0 )  |  
    w_1, w_2 \in {\bf R}_{\geq 0} \}.
$$
Here, the weight vector stands for
$(t_1,t_2,x,y, \partial_{t_1},\partial_{t_2}, \partial_x, \partial_y)$.

{\sl Gr\"obner fan in the homogenized Weyl algebra $\hooD$\/.}
We  consider the ideal $I_h$ in the homogenized Weyl algebra $\hooD$
generated by
\begin{eqnarray}  \label{eq:bs_example2}
&&h_\oo(g_1)=t_1-y, \ 
h_\oo(g_2)=t_2h-(yh-(x-h)^2), \\
&&h_\oo(g_3)=(-2x+2h)\partial_{t_2}+h \partial_{x}, \ 
h_\oo(g_4)=\partial_{t_1}+\partial_{t_2}+\partial_{y}. \nonumber 
\end{eqnarray}
Let $S$ be the stratification of $W$ by $\gr(\Dalg)$
or $\gr(\hzoDalg)$;
\begin{eqnarray*}
 S &=& \{ {\bar W'}, {\bar W}_{t_1}, {\bar W}_{t_2}, \{ 0 \} \}, \quad
     \mbox{ where } \\
 W'    &=& \{ (-w_1,-w_2, 0,0, w_1, w_2, 0,0 )  |  w_1, w_2 > 0 \} \\
 W_{t_1} &=& \{ (-w_1,0, 0,0, w_1, 0, 0,0 )  |  w_1 > 0 \} \\
 W_{t_2} &=& \{ (0,-w_2,0,0, 0, w_2, 0,0 )  |  w_2 > 0 \}. 
\end{eqnarray*}
The Gr\"obner fan ${\bar \E}(\hooD \cdot I_h, W) \cap S$ 
consists of 
$$ \{ {\bar F_1}, {\bar F_2}, {\bar L_{12}}, {\bar W_{t_1}},
      {\bar W_{t_2}}, \{ 0 \}  \} 
$$
where
\begin{eqnarray*}
  F_1 &=& \{ (-w_1,-w_2, 0,0, w_1, w_2, 0,0 ) \in W  |  
    w_1 >  w_2  \} \\
  F_2 &=& \{ (-w_1,-w_2, 0,0, w_1, w_2, 0,0 ) \in W  |  
    w_1 <  w_2  \} \\
  L_{12} &=& \{ (-w_1,-w_2, 0,0, w_1, w_2, 0,0 ) \in W  |  
    w_1 =  w_2  \}  
\end{eqnarray*}

The initial ideal for weights in $F_1$ is
$$
\ini_{F_1}(I_h) = 
 \{ -y, -x^2+2hx-h^2,
   2ht_2\partial_{t_2}+(hx-h^2)\partial_x +2h^3,
  (-2x+2h)\partial_{t_2}, \partial_{t_1} \}.
$$
The initial ideal for weights in $F_2$ is
$$
\ini_{F_2}(I_h) = 
 \{ -y,-x^2+2hx-h^2,
   2ht_1 \partial_{t_1}+(hx-h^2)\partial_x +2h^3,
   (-2x+2h)\partial_{t_1}, \partial_{t_2} \}.
$$
The initial ideal for weights in $L_{12}$ is
\begin{eqnarray*}
\ini_{L_{12}}(I_h)&=& 
 \{ 
  -y , 
  -x^2+2hx -h^2 , 
  2 h (t_1 \partial_{t_2}-t_2 \partial_{t_2})
    - (h x - h^2)\partial_x -2 h^3, \\
  && \quad \quad
  -2 (x \partial_{t_2}-h \partial_{t_2}), 
  \partial_{t_1}+\partial_{t_2} 
 \}.
\end{eqnarray*}
Here, $\ini_{F_1}(\cdot)$ means the initial with respect to a weight
vector in $F_1$.
Since the initial does not depend on the choice of the weight vector,
our notation does not have an ambiguity.

{\sl Gr\"obner fan in the Weyl algebra $D$\/.}
By utilizing the ideal membership algorithm, we can see that
$\ini_{F_1}(I_h)_{|h=1} \not= \ini_{F_2}(I_h)_{|h=1} 
 \not= \ini_{L_{12}}(I_h)_{|h=1}$.
Hence, we have
$$ {\bar \E}(D \cdot I, W)\cap S 
 = \{ {\bar F}_1,  {\bar F}_2, {\bar L}_{12},
      {\bar W}_{t_1}, {\bar W}_{t_2}, \{ 0 \} \}.
$$

{\sl Local Gr\"obner fan 
 in the local ring of differential operators $\Dalg$\/.}
First, let us compare 
$\ini_{F_1}(I_h)_{|h=1}$ and  $\ini_{F_2}(I_h)_{|h=1}$
in 
$\gr_{F_1}(\Dalg) = \{ \ell/g  |  \ell \in D, g \in {\bf Q}[x,y] \}$
By the tangent cone algorithm or by examining the bases, 
which contain the unit $-x^2+2x-1$, 
we can see that 
the dehomogenizations of the three initial ideals are 
the same ideal
$${\rm in}_{F_1} (I) = {\rm in}_{F_2}(I) = {\rm in}_{L_{12}}(I)
\quad {\rm in}\ \gr_{F_1} (\Dalg) = \gr_{F_2} (\Dalg) = \gr_{L_{12}}(\Dalg).
$$
The local Gr\"obner fan consists of only one maximal dimensional cone.
Hence, we have
$$ {\bar \E}(\Dalg \cdot I, W) = 
 \{ {\overline {F_1 \cup F_2}}, {\bar W}_{t_1}, {\bar W}_{t_2}, \{ 0 \} \} = S.
$$
This calculation shows us that 
either the global Gr\"obner fan 
${\bar \E} (D \cdot I, \Wloc) \cap {\bar S}_{\Dalg}$
or the local Gr\"obner fan 
${\bar \E} (\Dalg \cdot I, \Wloc)$ 
are not polyhedral fan.
Indeed, our comparison result Cor. \ref{cor:Wloc'} says that maximal
dimensional equivalence classes agree,
and if we suppose that they are both polyhedral fans,
then they must agree in the restricted weight space $W$, too.
It is a contradiction.

{\sl Global Gr\"obner fan in the homogenized Weyl algebra $\hzoD$\/.}
The algebra $\hzoD$ is the Weyl algebra defined with commutation
relation 
$ \pd_i x_j = x_j \pd_i + h \delta_{i,j} $.
Let us enumerate Gr\"obner cones of $I$ in the weight space $W$.
Note that $I$ is already $({\bf 0},{\bf 1})$-homogeneous.
We first homogenize $g_1, \ldots, g_4$ in the doubly homogenized
Weyl algebra defined with commutation relation
$ \pd_i x_j = x_j \pd_i + h h'\delta_{i,j} $ as follows;
\begin{eqnarray}  \label{eq:bs_example3} 
&&h'(g_1)=t_1-y, \ 
h'(g_2)=t_2 h'-(yh'-(x-h')^2), \\
&&h'(g_3)=(-2x+2h')\partial_{t_2}+h'\partial_{x}, \ 
h'(g_4)=\partial_{t_1}+\partial_{t_2}+\partial_{y}. \nonumber 
\end{eqnarray}
By a Gr\"obner basis computation of the ideal $J_1$
generated by the $h'(g_i)$'s, 
we get two maximal dimensional cones $F_1$ and $F_2$.
The initial ideals are
\begin{eqnarray*}
\ini_{F_1}(J_1) &=&
 \{    -y , \partial_{t_2} , -2 h' \partial_{t_1}+2 x \partial_{t_1} , 
    -{h'}^2+2 x h'-x^2 , \\
 & & \quad
  x h' \partial_x-2 t_1 x \partial_{t_1}-x^2 \partial_x-2 x h' h , 
  -4 t_1 x \partial_{t_1}^2-6 x^2 \partial_{t_1} h \},
\end{eqnarray*}
\begin{eqnarray*}
\ini_{F_2}(J_1) &=&
 \{ 
     -y , \partial_{t_1} , 2 h' \partial_{t_2}-2 x \partial_{t_2} , 
  -{h'}^2+2 x h'-x^2 ,  \\
 & & \quad
  x h' \partial_x-2 t_2 x \partial_{t_2}-x^2 \partial_x-2 x h' h , 
  -4 t_2 x \partial_{t_2}^2-6 x^2 \partial_{t_2} h  \},
\end{eqnarray*}
\begin{eqnarray*}
\ini_{L_{12}}(J_1) &=&
 \{ 
   -y , \partial_{t_1}+\partial_{t_2} , 2 h' \partial_{t_2}-2 x \partial_{t_2} , 
   -{h'}^2+2 x h'-x^2 ,  \\
& & \quad
xh'\partial_x+2 t_1 x \partial_{t_2}-2 t_2 x \partial_{t_2}-x^2 Dx-2 xh'h, \\
& & \quad
4 t_1 x \partial_{t_2}^2-4 t_2 x \partial_{t_2}^2-6 x^2 \partial_{t_2} h  
\}.
\end{eqnarray*}
By the ideal membership algorithm, we conclude that 
$\ini_{F_1}(I')_{|h'=1} \not= \ini_{F_2}(I')_{|h'=1}
 \not= \ini_{L_{12}}(I')_{|h'=1}$,
which implies that
$\ini_{F_1}(I) \not= \ini_{F_2}(I) \not= \ini_{L_{12}}(I)$
in $\hzoD$
and $F_1$, $F_2$ and $L_{12}$ are distinct.
Hence, we conclude that
$$ {\bar \E}(\hzoD \cdot I, W)\cap S 
 = \{ {\bar F}_1,  {\bar F}_2, {\bar L}_{12},
      {\bar W}_{t_1}, {\bar W}_{t_2}, \{ 0 \} \}.
$$

{\sl Local Gr\"obner fan in $\hzoDalg$\/.}
First, let us compare 
$\ini_{F_1}(I')_{|h=1}$ and  $\ini_{F_2}(I')_{|h=1}$
in 
$\gr_{F_1}(\hzoDalg) = 
 \{ \ell/g  |  \ell \in \hzoD, g \in {\bf Q}[x,y], g(0) \not= 0 \}$.
By the tangent cone algorithm or by examining the basis, 
we can see that 
the dehomogenizations of the two initial ideals are 
the same ideal in $\gr_{F_1}(\hzoDalg) = \gr_{F_2}(\hzoDalg)$.
The local Gr\"obner fan consists of only one maximal dimensional cone.
Hence, we have
$$ {\bar \E}(\hzoDalg \cdot I, W) = 
 \{ {\overline {F_1 \cup F_2}}, {\bar W}_{t_1}, {\bar W}_{t_2}, \{ 0 \} \} = S.
$$

This calculation shows us that 
the global Gr\"obner fan 
${\bar \E} (\hzoD \cdot I, \Wloc) \cap S$
is not a polyhedral fan.
Because, the local Gr\"obner fan 
${\bar \E} (\hzoDalg \cdot I, \Wloc)$ 
is a polyhedral fan and we have 
a comparison result Cor. \ref{cor:Wloc'Hom} which says that maximal
dimensional equivalence classes agree.

\end{ex}

\begin{rem*}\rm
In the computations above, the local Gr\"obner fans are trivial:
i.e. there is no slope between
$v_1=(-1, 0, 0, 0, 1, 0, 0, 0)$ and $v_2=(0, -1, 0, 0, 0, 1, 0, 0)$.
The global Gr\"obner fans are not trivial: there is
the slope $L_{12}=\langle (1,1) \rangle =
\langle 1\cdot v_1 + 1 \cdot v_2 \rangle $.
This may be explained by the following:
The only point where the local Bernstein-Sato
is not ``trivial'' (in some sense) is $(1,0)$
because for any $(x_0,y_0) \ne (1,0)$,
$f_1$ or $f_2$ becomes a unit in $\mathcal{O}_{(x_0, y_0)}$.
Now in local coordinates $(x',y')$ around $(1,0)$, $(f_1,f_2)$ is written
as $(y', y'-x'^2)$ and the Bernstein-Sato ideal is known to be
generated by $(s_1+1)(s_2+1)(2s_1 +2s_2+3)(2s_1+2s_2+5)$. We then
clearly see the linear form (or the slope) $(1,1)$.
\end{rem*}

\begin{ex} \rm
This example tries to show how far we can enumerate Gr\"obner cones
with our implementation.
Let $\ell_k$ be a differential operator
$$ \ell_k = \partial_{x_k} -
 \left( \sum_{i=1}^n x_i \partial_{x_i} + 1/2 \right)
 \left( x_k \partial_{x_k} + 1/(2k+1) \right).
$$
We consider the ideal 
$I^h$ in the homogenized Weyl algebra $\hooD$ in $n$ variables
generated by
$h_\oo(\ell_1), \ldots, h_\oo(\ell_n)$,
which is a system of differential equations for a Lauricella hypergeometric
series with irregular singularities when $h=1$. \\
\centerline{
\begin{tabular}{|l|l|}
\hline
$n$ & Number of maximal dimensional cones in ${\bf R}^{2n}$\\  \hline
$1$   & $2$  \\ \hline
$2$   & $39$  \\ \hline   
$3$   & $3246$  \\ \hline 
\end{tabular}
}
\end{ex}

\begin{ex} \rm
In our experience, computation in the doubly homogenized Weyl algebra 
$h'(D)$ sometimes exhausts huge memory space.
For example, consider the left ideal generated by
$$ x^3+y^2 h'-t_1 {h'}^2 , y^3+x^2 h'-t_2 {h'}^2 , 
 3 x^2 \pd{t_1}+2 x h' \pd{t_2}+{h'}^2 \pd{x} , 
 3 y^2 \pd{t_2}+2 y h' \pd{t_1}+{h'}^2 \pd{y}
$$ 
in the doubly homogenized Weyl algebra
${\bf Q}\langle h,h',x,y, t_1, t_2, \pd{x}, \pd{y}, \pd{t_1}, \pd{t_2} \rangle$.
We want to get the Gr\"obner fan in the restricted weight space
$$ \{ (0,0,-w_1,-w_2,0,0,w_1,w_2) \,|\, w_i \geq 0 \}, \quad
\mbox{where  $-w_i$ stands for $t_i$} 
$$
to compute a Bernstein-Sato polynomial for polynomials
$x^3+y^2$ and $x^2+y^3$ with the method given in \cite{Compos}.
However, our implementation exhausts 2G bytes of memory
in the stage of constructing the reduced Gr\"obner basis 
from a Gr\"obner basis of Collart-Kalkbrener-Mall's
Gr\"obner walk.
\end{ex}


\end{document}